\documentclass[draftclsnofoot,journal,onecolumn,12pt]{IEEEtran}
\usepackage{subfiles}
\usepackage{silence}
\usepackage{dsfont}

\usepackage[T1]{fontenc}
\usepackage{amsmath,amssymb,amsfonts,mathrsfs,bm}
\usepackage{mathtools}
\usepackage{amsthm}
\usepackage{cite}
\usepackage[shortlabels]{enumitem}
\usepackage{graphicx}
\usepackage{epstopdf}
\usepackage{url}
\usepackage{colortbl}
\usepackage{multirow}
\usepackage[table,dvipsnames]{xcolor}
\usepackage[normalem]{ulem}
\usepackage{array}
\newcolumntype{L}[1]{>{\raggedright\let\newline\\\arraybackslash\hspace{0pt}}m{#1}}
\newcolumntype{C}[1]{>{\centering\let\newline\\\arraybackslash\hspace{0pt}}m{#1}}
\newcolumntype{R}[1]{>{\raggedleft\let\newline\\\arraybackslash\hspace{0pt}}m{#1}}

\makeatletter
\let\MYcaption\@makecaption
\makeatother
\usepackage[font=footnotesize]{subcaption}
\makeatletter
\let\@makecaption\MYcaption
\makeatother

\usepackage{xparse}



\usepackage{glossaries}

\makeatletter
\let\oldgls\gls
\let\oldglspl\glspl

\newcommand\fussy@ifnextchar[3]{%
  \let\reserved@d=#1%
  \def\reserved@a{#2}%
  \def\reserved@b{#3}%
  \futurelet\@let@token\fussy@ifnch}
\def\fussy@ifnch{%
  \ifx\@let@token\reserved@d
    \let\reserved@c\reserved@a 
  \else
    \let\reserved@c\reserved@b
  \fi
 \reserved@c}

\renewcommand{\gls}[1]{%
  \oldgls{#1}\fussy@ifnextchar.{\@checkperiod}{\@}}
\renewcommand{\glspl}[1]{%
  \oldglspl{#1}\fussy@ifnextchar.{\@checkperiod}{\@}}

\newcommand{\@checkperiod}[1]{%
  \ifnum\sfcode`\.=\spacefactor\else#1\fi
}
\makeatother

\newacronym{wrt}{w.r.t.}{with respect to}
\newacronym{RHS}{RHS}{right-hand side}
\newacronym{LHS}{LHS}{left-hand side}
\newacronym{iid}{i.i.d.}{independent and identically distributed}

\usepackage{float}

\ifx\notloadhyperref\undefined
	\ifx\loadbibentry\undefined
		\usepackage[hidelinks,hypertexnames=false]{hyperref} 
	\else
		\usepackage{bibentry}
		\makeatletter\let\saved@bibitem\@bibitem\makeatother
		\usepackage[hidelinks,hypertexnames=false]{hyperref}
		\makeatletter\let\@bibitem\saved@bibitem\makeatother
	\fi
\else
	\relax
\fi

\usepackage[capitalize]{cleveref}
\crefname{equation}{}{}
\Crefname{equation}{}{}
\crefname{claim}{claim}{claims}
\crefname{step}{step}{steps}
\crefname{line}{line}{lines}
\crefname{condition}{condition}{conditions}
\crefname{dmath}{}{}
\crefname{dseries}{}{}
\crefname{dgroup}{}{}

\crefname{Problem}{Problem}{Problems}
\crefformat{Problem}{Problem~(#2#1#3)}
\crefrangeformat{Problem}{Problems~(#3#1#4) to~(#5#2#6)}

\crefname{Theorem}{Theorem}{Theorems}
\crefname{Corollary}{Corollary}{Corollaries}
\crefname{Proposition}{Proposition}{Propositions}
\crefname{Lemma}{Lemma}{Lemmas}
\crefname{Definition}{Definition}{Definitions}
\crefname{Example}{Example}{Examples}
\crefname{Assumption}{Assumption}{Assumptions}
\crefname{Remark}{Remark}{Remarks}
\crefname{Rem}{Remark}{Remarks}
\crefname{remarks}{Remarks}{Remarks}
\crefname{Exercise}{Exercise}{Exercises}
\crefname{Theorem_A}{Theorem}{Theorems}
\crefname{Corollary_A}{Corollary}{Corollaries}
\crefname{Proposition_A}{Proposition}{Propositions}
\crefname{Lemma_A}{Lemma}{Lemmas}
\crefname{Definition_A}{Definition}{Definitions}

\usepackage{algorithm,algpseudocode}

\ifx\loadbreqn\undefined
	\relax
\else
	\usepackage{breqn} 
\fi

\ifx\renewtheorem\undefined
\ifx\useTheoremCounter\undefined
\newtheorem{Theorem}{Theorem}
\newtheorem{Corollary}{Corollary}
\newtheorem{Proposition}{Proposition}
\newtheorem{Lemma}{Lemma}
\else
\newtheorem{Theorem}{Theorem}

\newtheorem{Proposition}[theorem]{Proposition}
\fi

\newtheorem{Definition}{Definition}

\fi

\theoremstyle{remark}

\theoremstyle{plain}

\newcommand{\Real}{\mathbb{R}}
\newcommand{\Nat}{\mathbb{N}}

\newcommand{\calA}{\mathcal{A}}

\newcommand{\calC}{\mathcal{C}}

\newcommand{\calG}{\mathcal{G}}

\newcommand{\calN}{\mathcal{N}}

\newcommand{\calR}{\mathcal{R}}

\newcommand{\calU}{\mathcal{U}}

\newcommand{\bbN}{\mathbb{N}}

\DeclareSymbolFont{bsfletters}{OT1}{cmss}{bx}{n}
\DeclareSymbolFont{ssfletters}{OT1}{cmss}{m}{n}
\DeclareMathSymbol{\bsfGamma}{0}{bsfletters}{'000}
\DeclareMathSymbol{\ssfGamma}{0}{ssfletters}{'000}
\DeclareMathSymbol{\bsfDelta}{0}{bsfletters}{'001}
\DeclareMathSymbol{\ssfDelta}{0}{ssfletters}{'001}
\DeclareMathSymbol{\bsfTheta}{0}{bsfletters}{'002}
\DeclareMathSymbol{\ssfTheta}{0}{ssfletters}{'002}
\DeclareMathSymbol{\bsfLambda}{0}{bsfletters}{'003}
\DeclareMathSymbol{\ssfLambda}{0}{ssfletters}{'003}
\DeclareMathSymbol{\bsfXi}{0}{bsfletters}{'004}
\DeclareMathSymbol{\ssfXi}{0}{ssfletters}{'004}
\DeclareMathSymbol{\bsfPi}{0}{bsfletters}{'005}
\DeclareMathSymbol{\ssfPi}{0}{ssfletters}{'005}
\DeclareMathSymbol{\bsfSigma}{0}{bsfletters}{'006}
\DeclareMathSymbol{\ssfSigma}{0}{ssfletters}{'006}
\DeclareMathSymbol{\bsfUpsilon}{0}{bsfletters}{'007}
\DeclareMathSymbol{\ssfUpsilon}{0}{ssfletters}{'007}
\DeclareMathSymbol{\bsfPhi}{0}{bsfletters}{'010}
\DeclareMathSymbol{\ssfPhi}{0}{ssfletters}{'010}
\DeclareMathSymbol{\bsfPsi}{0}{bsfletters}{'011}
\DeclareMathSymbol{\ssfPsi}{0}{ssfletters}{'011}
\DeclareMathSymbol{\bsfOmega}{0}{bsfletters}{'012}
\DeclareMathSymbol{\ssfOmega}{0}{ssfletters}{'012}

\DeclareMathOperator*{\argmax}{arg\,max}
\DeclareMathOperator*{\argmin}{arg\,min}

\DeclareMathOperator{\st}{s.t.}
\DeclareMathOperator{\as}{a.s.}

\DeclareMathOperator{\Tr}{Tr}

\DeclareMathOperator{\supp}{supp}

\DeclareMathOperator{\rank}{rank}

\DeclareMathOperator*{\esssup}{ess\,sup}

\DeclarePairedDelimiter\parens{(}{)}

\DeclarePairedDelimiter\braces{\{}{\}}

\newcommand{\qednew}{\nobreak \ifvmode \relax \else
      \ifdim\lastskip<1.5em \hskip-\lastskip
      \hskip1.5em plus0em minus0.5em \fi \nobreak
      \vrule height0.75em width0.5em depth0.25em\fi}

\newcommand{\nn}{\nonumber\\}

\newcommand{\indicator}[1]{{\bf 1}_{\braces*{#1}}}
\newcommand{\indicatore}[1]{{\bf 1}_{#1}}

\newcommand{\ofrac}[1]{{\frac{1}{#1}}}

\newcommand{\ceil}[1]{\left\lceil{#1}\right\rceil}

\newcommand{\KLD}[2]{{D({#1}\, \|\, {#2})}}

\DeclareDocumentCommand\set{s m t| m}{%
  \IfBooleanTF#1%
	{\left\{\, #2\mathrel{} \IfBooleanTF{#3}{\middle|}{:}\mathrel{}  #4\, \right\}}%
  {\{\, #2 \IfBooleanTF{#3}{\mid}{\mathrel{} : \mathrel{}} #4\, \}}%
}

\DeclareDocumentCommand \ifcond {m m} {%
	{#1} %
	\IfValueT{#2}{\, \middle|\, {#2}}%
}

\DeclareDocumentCommand \P {e{_} g >{\SplitArgument{ 1 }{ @| }}d() g } {%
	\mathbb{P}%
	\IfValueTF{#1}{_{#1}}
		{\IfValueT{#2}{_{#2}}}%
	\IfValueT{#3}{\left(\ifcond#3}%
	\IfValueT{#4}{\, \middle|\, {#4}}%
	\IfValueT{#3}{\right)}%
}

\DeclareDocumentCommand \E {e{_} g >{\SplitArgument{ 1 }{ @| }}o g } {%
	\mathbb{E}%
	\IfValueTF{#1}{_{#1}}
		{\IfValueT{#2}{_{#2}}}%
	\IfValueT{#3}{\left[\ifcond#3}%
	\IfValueT{#4}{\, \middle|\, {#4}}%
	\IfValueT{#3}{\right]}%
}

\definecolor{gray90}{gray}{0.9}

\ifx\nohighlights\undefined
	
	\newcommand{\blue}[1]{{{\color{blue} #1}}}
	\newcommand{\msout}[1]{\text{\color{green} \sout{\ensuremath{#1}}}}
	\newcommand{\del}[1]{{\color{green}\ifmmode \msout{#1}\else\sout{#1}\fi}}
\else
	
	\newcommand{\blue}[1]{#1}
	\newcommand{\msout}[1]{#1}
	\newcommand{\del}[1]{#1}
\fi

\newcommand{\hide}[1]{}

\renewcommand{\figurename}{Fig.}
\newcommand{\figref}[1]{\figurename~\ref{#1}}
\graphicspath{{./Figures/}} 
\ifx\diagnoselabel\undefined
	\relax
\else
	\makeatletter
	 \def\@testdef #1#2#3{%
		 \def\reserved@a{#3}\expandafter \ifx \csname #1@#2\endcsname
		\reserved@a  \else
	 \typeout{^^Jlabel #2 changed:^^J%
	 \meaning\reserved@a^^J%
	 \expandafter\meaning\csname #1@#2\endcsname^^J}%
	 \@tempswatrue \fi}
	\makeatother
\fi

\usepackage{algorithm}
\usepackage{algpseudocode}

\newcommand{\ARL}{\mathrm{ARL}}
\newcommand{\WADD}{\mathrm{WADD}}
\newcommand{\AOSC}{\mathrm{AOSC}}

\newcommand{\asure}[1]{\xrightarrow[\protect{\raisebox{3pt}[0pt][0pt]{\ensuremath{\scriptstyle\text{a.s.}}}}]{\P_{\nu,#1}}}
\newcommand{\asuren}{\xrightarrow[\protect{\raisebox{3pt}[0pt][0pt]{\ensuremath{\scriptstyle\text{a.s.}}}}]{\P}}

\newacronym{ARL}{ARL}{average run length}
\newacronym{WADD}{WADD}{worst-case average detection delay}
\newacronym{ADD}{ADD}{average detection delay}
\newacronym{AOSC}{AOSC}{average observation-switching cost}
\newacronym{QCD}{QCD}{quickest change detection}
\newacronym{QCQP}{QCQP}{quadratically constrained quadratic program}
\usepackage{xr}
\begin{document} 
	
\title{Asymptotically Optimal Sampling Policy for Quickest Change Detection with Observation-Switching Cost}

\author{
	Tze~Siong~Lau and Wee~Peng~Tay,~\IEEEmembership{Senior Member,~IEEE}
	\thanks{This work was supported in part by the Singapore Ministry of Education Academic Research Fund Tier 2 grant MOE2018-T2-2-019 and by A*STAR under its RIE2020 Advanced Manufacturing and Engineering (AME) Industry Alignment Fund – Pre Positioning (IAF-PP) (Grant No. A19D6a0053).}%
	\thanks{T.~S.~Lau and W.~P.~Tay are with the School of Electrical and Electronic Engineering, Nanyang Technological University, Singapore (e-mail: TLAU001@e.ntu.edu.sg, wptay@ntu.edu.sg).}
}

	\maketitle
	\begin{abstract}
		We consider the problem of quickest change detection (QCD) in a signal where its observations are obtained using a set of actions, and switching from one action to another comes with a cost. The objective is to design a stopping rule consisting of a sampling policy to determine the sequence of actions used to observe the signal and a stopping time to quickly detect for the change, subject to a constraint on the average observation-switching cost. We propose an open-loop sampling policy of finite window size and a generalized likelihood ratio (GLR) Cumulative Sum (CuSum) stopping time for the QCD problem. We show that the GLR CuSum stopping time is asymptotically optimal with a properly designed sampling policy and formulate the design of this sampling policy as a quadratic programming problem. We prove that it is sufficient to consider policies of window size not more than one when designing policies of finite window size and propose several algorithms that solve this optimization problem with theoretical guarantees. For observation-dependent policies, we propose a $2$-threshold stopping time and an observation-dependent sampling policy. We present a method to design the observation-dependent sampling policy based on open-loop sampling policies. Finally, we apply our approach to the problem of QCD of a partially observed graph signal and empirically demonstrate the performance of our proposed stopping times.
	\end{abstract}
	%

	\section{Introduction}\label{sec:intro}
	
	Quickest change detection (QCD) is the problem of detecting an abrupt change in a system while keeping the detection delay to a minimum. In a usual scenario, a sequence of \gls{iid} observations $\{x_t:t\in\mathbb{N}\}$, with probability density function (pdf) $p_0$ up to an unknown change point $\nu$, and i.i.d.\ with pdf $p_1\neq p_0$ after $\nu$, is obtained. The objective is to detect for the change at $\nu$ as quickly as possible while maintaining false alarm constraints \cite{tartakovsky2014sequential,poor2009quickest,veeravalli2013quickest}. QCD has applications across diverse fields, including quality control\cite{woodall2004using,lai1995sequential,lau17,LauTay:J19}, fraud detection\cite{bolton02}, cognitive radio\cite{lai2008quickest,ZhaTayLi:J16}, network surveillance\cite{tartakovsky2006novel,tartakovsky2014rapid}, structural health monitoring\cite{sohn00}, spam detection\cite{xie12,JiTayVar:J17,TanJiTay:J18}, bioinformatics\cite{muggeo10}, power system line outage detection\cite{banerjee2014power}, and sensor networks\cite{coppin1996digital,Hong2004,YanZhoTay:J18,lau19}.
	
	In many applications, the signal of interest $x_t$  may be high dimensional. For example, $x_t$ may consist of observations from many correlated sensors. Due to the large number of sensors in the network, bandwidth and power constraints prevent us from observing the entire network at any time instance, and we may only obtain sensor readings from a small subset of sensors at any time instance \cite{LiuTayLiu:J20,JiTay:J19}. While it may seem optimal to observe the maximum number of sensors allowed by the network, this sampling policy may not be feasible due to power and communication bandwidth considerations. Furthermore, the action of switching from one subset to another subset of sensors also incurs power and communication costs. In this paper, we consider both of these costs collectively as the \emph{observation-switching cost}, and we study the problem of QCD while maintaining an \gls{AOSC} constraint. To be more general, we consider the case where the signal can only be observed using an action selected from a set of permissible actions with observation-switching costs associated with the sequence of actions chosen. We assume that the pre- and post-change distributions as well as their conditional distributions given the actions are known to the observer. The objective is to design a sampling policy together with a stopping time that satisfies both the QCD false alarm and \gls{AOSC} requirements. To solve the QCD problem, we propose a sampling policy coupled with a generalized likelihood ratio (GLR) Cumulative Sum (CuSum) stopping time. For open-loop policies with finite window size, we show that the GLR CuSum stopping time is asymptotically optimal with a properly designed sampling policy and formulate the design of the sampling policy as a quadratic programming problem. For observation-dependent policies, we propose a $2$-threshold stopping time, prove that it satisfies the \gls{AOSC} and $\ARL$ constraints and demonstrate its performance empirically.
	
	\subsection{Related Work}
	Existing works in QCD where the signal is not entirely available to the decision maker or the fusion center can be categorized into three main categories. In the first category, the papers \cite{veeravalli1993decentralized,veeravalli2001decentralized,hadjiliadis2009one,jiang2016distributed} consider the problem of distributed or decentralized QCD where each node observes and processes its signal locally, with some memory of its previous messages, before sending a message to the fusion center. The authors of \cite{veeravalli2001decentralized} consider the problem where each sensor only has access to the local information at that node and would process the signal to send a quantized message to the fusion center for further processing.
	
	The second category of papers \cite{banerjee2012data,banerjee2013data,banerjee2013decentralized,geng2013non,banerjee15} consider the QCD problem where the number of observations made during the pre-change regime is controlled, and a control policy that determines whether a given observation is made. In \cite{banerjee15}, the authors developed a data-efficient scheme that allows for optional on-off sampling of the observations in the case where either the post-change family of distributions is finite, or both the pre- and post-change distribution belong to a one parameter exponential family.
	
	In the third category, the papers \cite{xie2015sketching,atia2015change,heydari2017} consider QCD where the observer only has access to compressed or incomplete measurements. The authors of \cite{xie2015sketching} study the problem of sequential change point detection where a randomly generated linear projection is used to reduce the dimensions of a high dimensional signal for the purpose of QCD. In \cite{heydari2017}, the authors consider QCD with an observation-dependent control of the actions where the current nodes to observe is determined by the maximal likelihood estimate of the post change hypothesis. In \cite{lau2017optimal}, we discussed the QCD problem where the observer is only able to obtain a partial observation of the signal through an action with an open-loop control of the actions.
	
	In the fourth category, the papers \cite{premkumar2008optimal,bajovic2011sensor,mei2011quickest,ren2016quickest} considers the QCD problem with observational scheduling considerations.  In these papers, it is assumed that there are multiple streams of observations and the cost associated with obtaining observations and the information quality of the streams differ from stream to stream. The QCD problem is to design a stopping time together with a control policy that determines the sequence of observations to perform such that the average cost of observations is controlled.
	
	Unlike the papers mentioned above, in this paper, we provide a general framework by considering random decision rules to select the current actions. We also do not give a fixed cost to the sampling of observation. Instead, we consider a more general approach where we use a set of permissible actions to model the practical sampling constraints and a cost is associated to the switching of actions to model observation-switching costs. In this paper, we consider the case where the decision maker is given a finite set of pre-defined actions, and the observed signal is a function of the action and the full signal. We also do not make any assumptions about the pre- and post-change distributions. 
	
	\subsection{Our Contributions}
	To the best of our knowledge, there are no existing work that considers the QCD problem while taking observation-switching costs into account. In this paper, we consider the problem QCD while maintaining an \gls{AOSC} constraint. The objective is to design a sampling policy together with a stopping time that satisfies both the QCD \gls{ARL} and \gls{AOSC} requirements. Our main contributions are as follows:
	\begin{enumerate}
		\item We formulate the QCD problem with an \gls{AOSC} constraint.
		\item We propose a condition on an open-loop policy of window size $W$ where the derivation of closed-form expressions for the \gls{AOSC} and the asymptotic \gls{WADD} of the GLR CuSum stopping time is possible.
		\item We prove the existence of an open-loop policy of window size $W$ satisfying our proposed condition for which the GLR CuSum is asymptotically optimal.
		\item Using these closed-form expressions, we formulate the design of the open-loop policy of window size $W$ as a quadratic programming problem with an additional combinatorial constraint and prove that for $1\leq W <\infty$, to obtain an asymptotically optimal policy, it is sufficient to consider only policies with window size 1.
		\item We propose an observation-dependent policy for which the GLR CuSum satisfies the \gls{AOSC} and \gls{ARL} requirements. We also experimentally verify the performance of the observation-dependent policy.
	\end{enumerate}
	
	The rest of this paper is organized as follows. In \cref{sec:problem}, we present our signal model and problem formulation. In \cref{sec:properties_AOSC_I}, we present properties of the \gls{AOSC} for open-loop sampling policies. In \cref{sec:static_policy_design}, we present the GLR CuSum stopping time and formulate the design of the open-loop sampling policy as a quadratic programming problem. We present algorithms to solve the open-loop policy design problems in \cref{sec:policy_design}. In \cref{sec:closed_loop_policy_design}, we present the $2$ threshold stopping time together with an observation-dependent sampling policy. Numerical results are presented in \cref{sec:results}. We conclude in \cref{sec:conclusion}.
	
	\emph{Notations:} The operator $\E_p$ denotes mathematical expectation \gls{wrt} pdf $p$, and $X \sim p$ means that the random variable $X$ has distribution with pdf $p$. The Kullback-Leibler (KL) divergence between the distributions with pdf $P$ and $Q$ is denoted as $\KLD{P}{Q}$. The pdf of a Gaussian distribution with mean $\mu$ and covariance $\Sigma$ is denoted as $\calN(\mu,\Sigma)$. Almost-sure convergence under the probability measure $\P$ is denoted as $\asuren$. We use $\indicatore{E}$ as the indicator function of the set $E$, and $\Nat$ to denote the set of positive integers. We use $\mathbb{R}$ to denote the set of real numbers and $\mathbb{R}^+$ to denote the set of positive real numbers. We also use the notation $a^{k:t}$ to denote the sequence $(a_k,a_{k+1},\ldots,a_t)$. For $\alpha=(a_1,a_2,\ldots,a_W)$, we use the notation $\alpha[j]=a_j$ to denote its $j$-th entry. For a probability transition matrix $T$, we use the notation $T[i,j]$ to denote the probability of moving to a state $j$ given that it is currently at state $i$. For a probability mass function $f$, $\supp(f)$ denotes the support of $f$.   
	
	\section{Problem formulation: Quickest Change Detection with a Cost for Switching Actions}\label{sec:problem}
	
	Let $p_0$ and $p_1,\ldots,p_M$ be $M+1$ distinct pdfs on $\Real^N$, and $X_1,X_2,\ldots$ be a sequence of vector-valued random variables satisfying the following:
	\begin{align}\label{eqn:signalmodel}
		\begin{cases}
			X_t \sim p_0  \text{ i.i.d.\ for all $t< \nu$},\\
			X_t \sim p_m \text{ i.i.d.\ for all $t\geq \nu$}.
		\end{cases}
	\end{align}
	where $\nu\geq 0$ and $m\in\{1,\ldots,M\}$ are unknown but deterministic constants. The QCD problem is to detect the change in distribution as quickly as possible by observing $X_1,X_2,\ldots$, while keeping the false alarm rate low. 
	
	In this paper, we assume that the observer is only able to obtain a partial observation $(A_t,Y_t)$ of $X_t$, where $Y_t\triangleq A_t(X_t)$ is a function of the random variable $X_t$ under the action $A_t$ at each time $t$. Let $\mathcal{A}$ be the collection of permissible actions. We assume that the set $\mathcal{A}=\{1,2,\ldots,|\calA|\}$ is finite. We also assume that at each time $t$, under the pdf $p_m$, the observation $Y_t$ is conditionally independent of $Y_1,\ldots,Y_{t-1}$ and $A_1,\ldots,A_{t-1}$ given the action $A_t$. We further assume that for each $m\in\{1,\ldots,M\}$ there exists an action $A\in\calA$ such that $p_0$ and $p_m$ are distinguishable under the action $A$. Some examples of $\mathcal{A}$ that arise in practical applications include:
	\begin{enumerate}
		\item 
		Network Sampling. The set of rank $n<N$ transformations with
			\begin{align*}
				&\mathcal{A} = \left\{A_{\mathbf{L}}:\mathbb{R}^N\to\mathbb{R}^n \ \middle\vert \begin{array}{l}
					A_{\mathbf{L}}(X)=\mathbf{L}X,\\
					\mathbf{L}=[e_{i_1},\ldots,e_{i_n}]^T,\\
					i_1<i_2<\ldots<i_n
				\end{array}\right\}
			\end{align*}
			where $e_i$ is an $N\times 1$ column vector with all zeros except a $1$ at the $i$-th position. 
		\item 
		$1$-bit quantization. A set of $2^n$ functions on $\mathbb{R}^+$ with 
			\begin{align*}
				&\mathcal{A}=\left\{A_{\phi}:\mathbb{R}^+\to \{0,1\}\ \middle\vert \begin{array}{l}
					A_{\phi}(x)=\indicator{x\geq \phi},\\
					\phi\in\{1,2,3,\ldots,2^n\} 
				\end{array}\right\}
			\end{align*}
			where $n$ is a positive integer, $\indicatore{A}$ is the indicator function of the set $A$ and $\phi$ is the quantization threshold.
	\end{enumerate}
	
	In our sequential change detection problem, we obtain observations $(A_1=a_1,Y_1=y_1),(A_2=a_2,Y_2=y_2),\ldots$ sequentially and aim to detect the change in pdf from $p_0$ to $p_m$ for some fixed $m\in\{1,...,M\}$ as quickly as possible. This is determined by a stopping time. A sampling policy is used to determine the action used to obtain the next observation. An observation-dependent sampling policy determines the current action based on both the previous actions and observations while an open-loop sampling policy determines the current action based on only the previous actions. Since an observation-dependent sampling policy has access to more information, an optimized observation-dependent policy is expected to outperform an open-loop policy. 
		
		 However, an observation-dependent policy may not be suitable for some applications due to physical constraints. One such class of applications involves time-synchronized sensor networks with no feedback mechanism from the fusion center to the sensors. In this case, each sensor in a time-synchronized network may have several modes (temperature, humidity, pressure, etc.) of obtaining observations. Without a feedback mechanism from the fusion center, the only way for the sensor to select a mode of observation is to use an open-loop policy. 
		Another class of applications is when the sensor network has a transmission delay that is significantly larger the sampling rate of the sensors. When an observation-dependent policy is applied, the sampling rate of this network is greatly 	reduced due to the time taken to communicate actions to the sensors. An open-loop policy does not suffer from the same problem as the sequence of actions does not depend on the observations, and thus can be pre-generated and made available to the sensors. Geo-stationary satellites usually have a communication round-trip time of between $600-800$ms while carrying sensors capable of sampling at rates larger than 1kHz. An observation-dependent policy would only be able to obtain a sample from the geo-stationary satellites at sampling frequency of less than 2Hz since it has to wait for the observations to be transmitted from the satellites to the ground station before deciding on the next action. An open-loop policy does not need to wait for the observations to decide the next action as the policy has been pre-determined. Hence, an open-loop policy would be able to sample the network of satellite at rates larger than 1kHz with a delay of between $300-400$ms from the time the observation is obtained to the time the observation is processed by the stopping time. In this example, it is possible that a stopping time using an open-loop policy out-performs one that uses an observation-dependent policy due to large differences in the amount of observations available. 
		In both classes, a long sequence of actions may be pre-generated and programmed into the sensors and the fusion center so that there is no need to communicate the sequence of actions during test time. In this paper, we consider only open loop-sampling policies. The interested reader may refer to the supplementary material \cite{lau2020supplementary} for a discussion on observation-dependent sampling policies. 
	\begin{Definition}
		A \emph{policy} $\pi$ is a sequence of functions $(\rho_t)_{t\in\mathbb{N}}$, where $\rho_t$ is a randomized function that determines the action $A_{t+1}$ using observations $\left((A_1,Y_1),(A_2,Y_2),\ldots,(A_t,Y_t)\right)$ up to time $t$.\\
		An \emph{open-loop policy $\pi$ of window size $W$} is a policy $(\rho_t)_{t\in\mathbb{N},t \geq W}$, where $\rho_t=\rho$ for a fixed randomized function $\rho$ that determines the action $A_{t+1}$ based on $W$ past actions $\left(A_{t-W+1},A_{t-W+2},\ldots,A_t\right)$.
	\end{Definition}
	It can be shown that an open-loop policy $\pi$ of window size $W$ is equivalent to a Markov chain of order $W$ on $\calA$ with initial distribution $q$ and probability transition matrix $T$. Thus, we use the notation $\pi=(q,T)$ to represent an open-loop policy of window size $W$ for the rest of this paper. If the change point is at $\nu$ and post-change distribution has pdf $p_m$, we let $\P{\nu,m}$ and $\E_{\nu,m}$ be the probability measure and mathematical expectation, respectively. We let $\P{\infty}$ and $\E_{\infty}$ denote the probability measure and mathematical expectation when there is no change.

	
	For a stopping time $\tau$ and a policy $\pi$, we quantify its detection delay using the worst case average detection delay (WADD) as proposed by Lorden\cite{lorden71}:
	\begin{align}
		&\text{WADD}(\tau,\pi)=\max_{1\leq m\leq M}\text{WADD}_m(\tau,\pi),
	\end{align}
	where
	\begin{align}
		&\text{WADD}_m(\tau,\pi) \nonumber\\
		&=\sup_{\nu\geq 1}\esssup \E{\nu,m}[(\tau-\nu+1)^+]{A^{1:(\nu-1)},Y^{1:(\nu-1)}},
	\end{align} 
	and its \gls{ARL} to false alarm as $\text{ARL}(\tau,\pi)=\E_{\infty}[\tau]$.
	
	In order to take the observation-switching costs of a policy $\pi$ into consideration, we let $\calC$ be a $|\calA|\times|\calA|$ matrix where its $(i,j)$-th entry $\calC[i,j]$ denotes the cost of switching from action $i$ to action $j$. Inspired by a similar cost first proposed in \cite{banerjee2013data}, we define the \gls{AOSC} of the policy $\pi$ as
	\begin{align}\label{eqn:def:AOSC}
		\text{AOSC}(\pi)=\limsup_{n\to \infty}\frac{1}{n}\E{\infty}[\sum_{t=2}^{n+1} \calC[A_{t-1},A_{t}]].
	\end{align} 
	Formally, our quickest change detection with \gls{AOSC} constraint can be formulated as a minimax problem: find a sampling policy $\pi$ and a stopping time $\tau$ to
	\begin{align}\label[Problem]{eqn:QCD_AOSC_formulation}
		\begin{aligned}
			& \min_{\tau,\pi} &  & \WADD(\tau,\pi)           \\
			& \st       &  & \ARL(\tau,\pi)\geq \gamma, \\
			&      &  & \AOSC(\pi)\leq\alpha_{\AOSC},
		\end{aligned}
	\end{align} 
	for some given thresholds $\alpha_{\text{AOSC}}$ and $\gamma$. 
	
	For a fixed policy $\pi$, the GLR CuSum\cite{banerjee15} stopping time $\tau_\pi$ \gls{wrt} the observed sequence $(A_t=a_t,Y_t=y_t)_{t\in\bbN}$ is defined as: 
	\begin{align}
		\tau_\pi &= \inf\set* t | {S(t,\pi)>\log (M\gamma)},\\
		S(t,\pi)&=\max_{1\leq m\leq M} \max_{1\leq i \leq t+1}\sum_{j=i}^{t}\log \frac{p_m(y_j|a_j)}{p_0(y_j|a_j)},  \label{eqn:parallel_CuSum_stoptime}
	\end{align}
	where $\gamma\geq0$ is a pre-selected threshold and $p_m(y_j|a_j)$ is the conditional pdf of $Y_j=y_j$ given the action $A_j=a_j$ under the distribution with pdf $p_m$. The GLR CuSum stopping time $\tau_\pi$ can be re-written as  
	\begin{align}
		\tau_\pi&=\min_{1\leq m\leq M}\tau_{m,\pi},\quad S(t,\pi)=\max_{1\leq m\leq M} S_m(t,\pi),\\
		\tau_{m,\pi}&=\inf\set*{t}|{S_m(t,\pi)>\log (M \gamma)},\\
		S_m(t,\pi)&=\left(S_m(t-1,\pi)+\log \frac{p_m(y_t|a_t)}{p_0(y_t|a_t)}\right)^+\text{for $t>0$},\label{eqn:parallel_CuSum_update}\\
		S_m(0,\pi)&=0,\quad\text{for $m\in\{1,\ldots,M\}$}
	\end{align}
	where $x^+\triangleq\max(x,0)$. We note that for $m\in\{1,\ldots,M\}$, $S_m(t,\pi)$ is the CuSum statistic corresponding to the post-change pdf $p_m$ and policy $\pi$. Thus, the GLR CuSum statistic $S(t,\pi)$ is the maximum of the CuSum statistics for each of the post-change pdf $p_m$.
	
	We discuss results pertaining to open-loop policies of window size $W$ in \cref{sec:policy_design,sec:properties_AOSC_I,sec:static_policy_design}, and propose an observation-dependent policy in \cref{sec:closed_loop_policy_design}.
	
	\section{Properties of the AOSC for open-loop policies of window size $W$}\label{sec:properties_AOSC_I}
	In this section, we present results regarding the $\AOSC$ of an open-loop policy of window size $W$. When the window size $W$ of the open loop policy is zero, the actions used to observe the signal are generated \gls{iid} with respect to the distribution $q$. In this case, the observations  $\{(A_t,Y_t)=(a_t,y_t)\}_{t\in\mathbb{N}}$ are also generated \gls{iid}. However, when the window size $W$ of the open-loop policy is positive, unlike the former case, it is possible that the actions $\{A_t\}_{t\in\mathbb{N}}$ and observations $\{(A_t,Y_t)\}_{t\in\mathbb{N}}$ are not generated \gls{iid}. We denote the joint probability density function of $(A^{1:t},Y^{1:t})$ under the distribution with pdf $p_m$ as 
	\begin{align}
		&p_m(a^{1:t},y^{1:t})\nonumber\\
		&= q(a^{1:W})\prod_{j=1}^{t}p_m(y_j|a_j)\prod_{k=W+1 }^{t}p_{T}(a_k\ |\ a^{k-W:k-1}),
	\end{align}
	where $p_{T}$ is the conditional probability mass function of $A_k$ given $A^{k-W:k-1}$ induced by the probability transition matrix $T$, $W$ is the window size of the policy.
	
	An open-loop policy $\pi=(q,T)$ with window size $W$ can also be written as a Markov chain $\pi'=(q',T')$ of order $1$ where $T'$ satisfies $T'[\alpha,\beta]=0$ whenever $\beta[i]\neq \alpha[i-1]$ for some $i\in \{2,\ldots,W\}$ and $\alpha,\beta\in\calA^W$. For the rest of this paper, we switch between either representation of an open-loop policy $\pi$ to simplify the computations in the proofs. We denote the observation-switching costs associated with the latest two actions $\calC[\alpha[W-1],\alpha[W]]$, from $\alpha[W-1]$ to $\alpha[W]$, as $\calC_\alpha$.
	
	For the rest of this section, we present results that relate the $\AOSC$ of open-loop polices with different initial distributions but equal probability transition matrices. First, we recall a relation between the average number of visits and the stationary distributions of a Markov chain.  Let $N_t(\alpha;\beta)$ denote the number of times, up to time $t$, that the state $\alpha$ is visited given that the initial state is $\beta$. Since $\calA^W$ is finite, the Markov chain defined by the transition matrix $T$ has at least one recurrence class. Let $R$ be the number of recurrence classes and $U$ be the number of transient states. By the Ergodic Theorem for finite state Markov chains\cite{bertsekas2002introduction}, for a finite state Markov Chain with $R$ recurrent classes $\{\calR_1,\calR_2,\ldots,\calR_R\}$, there exists $R$ stationary distributions $\xi_1,\ldots\xi_R$ where $\xi_r[\alpha]=0$ if the state $\alpha\notin \calR_r$, and for recurrent states $\beta\in \calR_r$, we have 
	\begin{align}
		&\lim_{t\to\infty}\frac{N_t(\alpha;\beta)}{t} = \xi_r[\alpha]\ \as, \\
		& \lim_{t\to\infty}\frac{\E[N_t(\alpha;\beta)]}{t} = \xi_r[\alpha],
	\end{align}
	for any state $\alpha\in\calA^W$ and $r=1,\ldots,R$. For transient states $\beta$, 
	\begin{align}
		\lim_{t\to\infty}\frac{\E[N_t(\alpha;\beta)]}{t}=\sum_{r=1}^R f_{\beta,r}\xi_r[\alpha],
	\end{align}
	where $f_{\beta,r}$ is the first-passage probability of initializing at state $\beta$ and entering the recurrence class $\calR_r$ before any other recurrence classes.
	
	For any state $\beta$, denoting $\xi$ as the vector of expected proportion of visits to each of the states initializing at state $\beta$ such that $\xi[\alpha]=\lim_{t\to\infty}\frac{\E[N_t(\alpha;\beta)]}{t} $, we can see that $\xi$ is a stationary distribution of the probability transition matrix $T$ as it is a convex linear combination of stationary distributions. 
	\begin{Definition}
		For any initial distribution $q$, the expected proportion of visits to each of the state, $\overline{q}$, is defined as 
		\begin{align*}
		\overline{q}(\alpha)=\lim_{t\to\infty}\E\left[\frac{1}{t}\sum_{j=1}^t\mathbf{1}_{\{A_j=\alpha\}}\right].
		\end{align*}
	\end{Definition}
Thus, for any initial distribution $q$, $\overline{q}$ is a stationary distribution of $T$. In the next lemmas, we see that the AOSC and asymptotic log likelihood ratios depend only on $T$ and $\overline{q}$.
	
	\begin{Lemma}\label{lem:AOSC_staionary}
		Let $\pi_1=(q,T)$ be an open-loop policy of finite window size $W$ and $\pi_2=(\overline{q},T)$, then we have 
		\begin{align}
			\AOSC(\pi_1)=\AOSC(\pi_2)
		\end{align}
	\end{Lemma} 
	\begin{IEEEproof}
		The proof is based on standard Markov chain theory and is provided for completeness in the supplementary material \cite{lau2020supplementary}.
	\end{IEEEproof}

	\section{Asymptotic Properties of GLR CuSum for open-loop policies of window size $W$}\label{sec:static_policy_design}
	
	Next, we present some asymptotic properties of $S(t,\pi)$ and $\tau_{\pi}$ for an open-loop policy $\pi$. In this paper, we use $\asymp$ to denote the notion of asymptotic equivalence\cite{greene2007mathematics}:
	\begin{align}
		f\asymp g \quad \text{if and only if}\quad  \lim_{x\to \infty} \frac{f(x)}{g(x)}=1.
	\end{align}
		\begin{Definition}\label{def:trade-off_rate}
			The \emph{asymptotic $\ARL$-$\WADD$ trade-off rate} $I$ of the GLR CuSum stopping time $\tau_{\pi}$ is defined as 
			\begin{align*}I=\liminf_{\gamma\to\infty}\frac{\log \ARL(\tau_{\pi}(\gamma))}{\WADD(\tau_{\pi}(\gamma),\pi)}.\end{align*}
		\end{Definition}
		When the signal is generated i.i.d. before and after the change point, the asymptotic trade-off rate for the GLR CuSum stopping time is well studied and has a nice closed-form expression in terms of the KL divergence between the pre-change and post-change distributions \cite{banerjee15}.
		In our case, even though the signal $X_1,X_2,\ldots$ is originally i.i.d. before and after the change-point, any sampling procedure that switches action would inevitably result in non-i.i.d. observations.
	    We let 
		\begin{align}
			\Lambda_m(k,t)=\log\frac{p_m(A^{k:t},Y^{k:t})}{p_0(A^{k:t},Y^{k:t})}=\sum_{i=k}^t\log\frac{p_m(Y_i|A_i)}{p_0(Y_i|A_i)}.
		\end{align}
		In the next two Lemmas, we present results regarding $\Lambda_m$ when the open-loop policy $\pi=(q,T)$  satisfies the property that $\overline{q}$  has support in only one recurrence class of $T$. The property that $\overline{q}$  has support in only one recurrence class of $T$ plays an important role in solving \cref{eqn:QCD_AOSC_formulation}. It allows us to quantify the performance of the GLR CuSum using a closed form expression. This is key to obtaining asymptotically optimal policies using numerical optimization tools.
		
		\begin{Lemma}\label{prop:information_number_single_recurrence_as_convergence}
			For any open-loop policy $\pi=(q,T)$ of finite window size $W$ where $\overline{q}$ has support in only one recurrence class $\calR$, and any $m \in\{1,\ldots,M\}$ and change-point $\nu<\infty$, we have 
			\begin{align}\label{eqn:information_number_as_convergence}
				\P{\nu,m}(\lim_{t\to\infty}\frac{1}{t}\Lambda_m(\nu,\nu+t-1)= I_{m,\pi})=1,
			\end{align}
			with 
			\begin{align}\label{Impi}
				I_{m,\pi}=\sum_{\alpha\in \calR}\overline{q}(\alpha) \KLD{p_m(\,\cdot\mid\alpha[W])}{p_0(\,\cdot\mid\alpha[W])}.
			\end{align}
		\end{Lemma}
		\begin{IEEEproof}
			See \cref{sec:AppLem2}.
		\end{IEEEproof}
		
		\begin{Lemma}\label{prop:information_number_single_recurrence_lower_bounds_LLR}
			Let $\pi=(q,T)$ be an open-loop policy of finite window size $W$ where $\overline{q}$ has support in only one recurrence class. For any $\epsilon>0$ and $m \in\{1,\ldots,M\}$, we have
			\begin{align}
				&\lim_{t\to\infty}\P{\nu,m}(\left|\frac{1}{t}\Lambda_m(\nu,\nu+t-1)- I_{m,\pi}\right|>\epsilon)=0, \label{eqn:information_number_convergence_in_prob}
			\end{align}
			for $0\leq \nu <\infty$, and
			\begin{multline}
				\sup_{0\leq \nu<\infty}\esssup \mathbb{P}_{\nu,m}\bigg(\ofrac{t}\max_{0\leq j<t} \Lambda_m(\nu,\nu+j)\\
				>(1+\epsilon)I_{m,\pi}\ \bigg| \ {A^{1:\nu-1},Y^{1:\nu-1}}\bigg)\to 0\quad\text{as $t\to\infty$.}\label{eqn:information_number_convergence_in_prob_lower_bound}
			\end{multline}
		\end{Lemma}
		\begin{IEEEproof}
			See \cref{sec:AppLem3}.
		\end{IEEEproof}
		
	 For a fixed open-loop policy $\pi=(q,T)$ such that the expected proportion of visits to each of the states, $\overline{q}$, has support in one recurrence class, we apply \cref{prop:information_number_single_recurrence_lower_bounds_LLR} together with \cite[Theorem 8.2.3]{tartakovsky2014sequential} to obtain the following proposition.
	
	\begin{Proposition}\label{prop:Information_number_WADD}
		For a fixed open-loop policy $\pi=(q,T)$ such that $\overline{q}$ has support in one recurrence class, we have the following asymptotic $\ARL$-$\WADD$ trade-off for any $1\leq m \leq M$:
		\begin{align}
			\ARL(\tau_{m,\pi}, \pi)\geq M\gamma,\ 
			\WADD_m(\tau_{m,\pi}, \pi)\leq \frac{\log \gamma}{I_{m,\pi}}(1+o(1)),
		\end{align} 
		as $\gamma\to\infty$ where a function $g(\gamma)=o(1)$ if and only if $\lim_{\gamma\to\infty}g(\gamma)=0$.
	\end{Proposition}
	Thus, when the signal $\{X_t\ :\ t\in\mathbb{N}\}$ is sampled using the open-loop policy $\pi$ such that $\overline{q}$ has support in one recurrence class, using similar techniques from \cite[Theorem 6.16]{poor2009quickest} together with \cref{prop:Information_number_WADD}, we know that the GLR CuSum stopping time $\tau_{\pi}$ gives us a stopping time satisfying $\ARL(\tau_{\pi},\pi)\geq \gamma$ and $\tau_{\pi}$ is asymptotically optimal for the following problem as $\gamma\to\infty$: 
	\begin{equation}\label{eqn:QCD_formulation_specific_sampling}
	\begin{aligned}
	& \min_{\tau} &  & \WADD(\tau, \pi)& \st   &  & \ARL(\tau, \pi)\geq \gamma.
	\end{aligned}
	\end{equation}
	It should be noted that, for open-loop policies $\pi=(q,T)$ with $\overline{q}$ having support in more than one recurrence class, \cref{prop:Information_number_WADD} does not hold, making it difficult to use \cite[Theorem 6.16]{poor2009quickest} to derive the asymptotic $\ARL$-$\WADD$ trade-off rate and show the asymptotic optimality of the GLR CuSum. Furthermore, in the next proposition, we show that an open-loop policy $\pi$ with $\overline{q}$ having support in multiple recurrence classes is suboptimal in terms of $\AOSC$ and $\WADD$.
	
	\begin{Proposition}\label{prop:reduction_of_recurrence_classes}
		Let the open-loop policy $\pi=(q,T)$ be such that $\overline{q}$ has support in multiple recurrence classes. Then, there exists an open-loop policy $\pi'=(q',T)$ where $\overline{q'}$ has support in only one recurrence class such that for any stopping time $\tau$,
		\begin{align*}
			\AOSC(\pi')\leq \AOSC(\pi) \text{ and } \WADD(\tau,\pi')\leq \WADD(\tau,\pi).
		\end{align*}
	\end{Proposition}
	\begin{IEEEproof}
		See \cref{sec:AppProp2}.
	\end{IEEEproof}
	
	Using this proposition, we obtain a result regarding asymptotically optimal solutions of \cref{eqn:QCD_AOSC_formulation}.
	\begin{Theorem}\label{thm:structure_of_asym_opt_solutions}
		When the signal is sampled using an open-loop policy $\pi=(q,T)$ with $q$ having support in one recurrence class, satisfying $\AOSC(\pi)\leq\alpha_{\text{AOSC}}$, the GLR CuSum $\tau_{\pi}$ is asymptotically optimal with the asymptotic $\WADD$-$\ARL$ trade-off given as $\min_{m} I_{m,\pi}$. 
			
		There exists an open-loop policy $\pi=(q',T')$ with $\overline{q'}$ having support in one recurrence class such that $(\tau_\pi,\pi)$ is asymptotically optimal for \cref{eqn:QCD_AOSC_formulation} as $\gamma\to \infty$, where $\tau_\pi$ is the GLR CuSum stopping time. 
	\end{Theorem}
	\begin{IEEEproof}
		See \cref{sec:AppThm1}.
	\end{IEEEproof}
	
	Thus, the additional constraint that the open-loop policy $\pi=(q,T)$ satisfies the condition that $\overline{q}$ has support in one recurrence class does not affect the asymptotic optimality of the GLR CuSum for \cref{eqn:QCD_AOSC_formulation}. Furthermore, for policies satisfying this constraint, we are able to obtain a closed-form expression for the $\ARL$-$\WADD$ trade-off rate. The closed-form expressions are important as numerical optimization tools are used to optimize the $\ARL$-$\WADD$ trade-off for the GLR CuSum.
	
	Using \cref{thm:structure_of_asym_opt_solutions}, the minimization in \cref{eqn:QCD_AOSC_formulation}, over the open-loop policy $\pi$ and stopping time $\tau$, can be decoupled.
	
	Let $\pi^*$ be an optimal solution to the following problem:
	\begin{align}\label[Problem]{eqn:optimize_sampling_policy_asymptotic}
		\begin{aligned}
			& \min_\pi & & \max_{1\leq m \leq M}\ I_{m,\pi}^{-1}\\
			&\st && \AOSC(\pi)\leq\alpha_{\AOSC},\\
			& &&\supp(\overline{q}) \subseteq \text{one recurrence class}.
		\end{aligned}
	\end{align}
	By the argument above, $(\tau_{\pi^*},\pi^*)$ is asymptotically optimal for \cref{eqn:QCD_AOSC_formulation}. 	We call \cref{eqn:optimize_sampling_policy_asymptotic} the \emph{open-loop policy design} problem.

	\section{Optimal open-loop policy of window size $W$}\label{sec:policy_design}
	
	In this section, we investigate the open-loop policy design \cref{eqn:optimize_sampling_policy_asymptotic} under the cases where the switching costs from one action to another are all equal or not.
	
	\subsection{Equal Switching Costs}\label{subsect:policy_design_equal_cost}
	In this subsection, we propose a method to solve the open-loop policy design problem in which the switching costs are constant, i.e., $\calC[a,b]=c$, for a fixed $c\in\mathbb{R}$ and any $a,b\in \calA$. First, we note that \cref{eqn:optimize_sampling_policy_asymptotic} is feasible if and only if $\alpha_{\text{AOSC}}\geq c$. Furthermore, if $\alpha_{\text{AOSC}}\geq c$ then \cref{eqn:QCD_AOSC_formulation} reduces to 
	\begin{equation}\label[Problem]{eqn:QCD_AOSC_formulation_remove_AOSC}
	\begin{aligned}
	& \min_{\tau,\pi} &  & \WADD(\tau,\pi)           \\
	& \st       &  & \ARL(\tau,\pi)\geq \gamma,
	\end{aligned}
	\end{equation}
	where the $\AOSC$ constraint is automatically satisfied. Next, we show that for the case when all action-switching costs are equal, there exists a memoryless open-loop policy $\pi$ (i.e., $W=0$) for which the GLR CuSum $\tau_{\pi}$ achieves asymptotic optimality.
	\begin{Proposition}\label{prop:M=0_suffices}
		Suppose \cref{eqn:optimize_sampling_policy_asymptotic} is feasible and $(\tau_{\pi},\pi)$ is an asymptotically optimal solution for \cref{eqn:QCD_AOSC_formulation}. Then, there exists an open-loop policy $\pi_0$ with window size $W=0$ such that
		\begin{align}
			\WADD(\tau_{\pi_0},\pi_0)\asymp\WADD(\tau_{\pi},\pi)\quad\text{as $\gamma\to\infty$.}
		\end{align}
	\end{Proposition}
	\begin{IEEEproof}
		See \cref{sec:AppProp4}.
	\end{IEEEproof}
	
	From \cref{prop:M=0_suffices}, to solve the open-loop policy design problem \cref{eqn:optimize_sampling_policy_asymptotic} for some $W\in\mathbb{N}$, it suffices to solve \cref{eqn:optimize_sampling_policy_asymptotic} for the case where $W=0$. 
	
	When $W=0$, \cref{eqn:optimize_sampling_policy_asymptotic} becomes
	\begin{equation}\label[Problem]{eqn:policy_design_problem_m=0}
	\begin{aligned}
	& \min_q
	& & \max_{1\leq m\leq M} \parens*{\sum_{a\in \mathcal{A}} q(a) \KLD{p_0(\,\cdot\mid a)}{p_m(\,\cdot\mid a)}}^{-1}\\
	& \st
	& & \sum_{a\in\mathcal{A}}q(a)=1, \\
	& & & q(a)\geq 0 \quad \text{for all $a\in\mathcal{A}$}.
	\end{aligned}
	\end{equation}
	
		\begin{Proposition}\label{prop:equivalence_of_optimization_problem}
			For the optimization problem
			\begin{equation}\label[Problem]{eqn:policy_design_problem_m=0_LP}
			\begin{aligned}
			& \min_{q,z} & &z\\
			& \st
			& & \sum_{a\in\mathcal{A}}q(a)=1,\\
			& & & q(a)>0 \quad \text{for all $a\in\mathcal{A}$},\\
			& & &\sum_{a\in \mathcal{A}} q(a) \KLD{p_0(\,\cdot\mid a)}{p_m(\,\cdot\mid a)}+z\geq 0\\
			& & &\quad\quad\text{for all $m\in\{1,\ldots,M\}$},
			\end{aligned}
			\end{equation} the following holds
			\begin{enumerate}
				\item[(i)] If $(q^*,z^*)$ is an optimal solution to \cref{eqn:policy_design_problem_m=0_LP} then $q^*$ is an optimal solution to \cref{eqn:policy_design_problem_m=0}.
				\item[(ii)] If  $q^*$ is an optimal solution to \cref{eqn:policy_design_problem_m=0} then there exists $x^*$ such that $(q^*,x^*)$ is an optimal solution to \cref{eqn:policy_design_problem_m=0_LP}.
			\end{enumerate}
		\end{Proposition}
		
		\begin{IEEEproof}
			See \cref{sec:AppProp4.5}.
		\end{IEEEproof}
	
	By \cref{prop:equivalence_of_optimization_problem}, we can solve \cref{eqn:policy_design_problem_m=0} by solving the linear optimization problem, \cref{eqn:policy_design_problem_m=0_LP}. Let $q_0$ be the solution for Problem \eqref{eqn:policy_design_problem_m=0_LP} and $T_0$ be the probability transition matrix with rows equal to $q_0$. Using similar techniques from \cite[Theorem 6.16]{poor2009quickest} together with \cref{prop:Information_number_WADD}, we know that the GLR CuSum algorithm with $\pi_0=(q_0,T_0)$ as the open-loop policy gives us a stopping time satisfying $\text{ARL}(\tau_{q^*})\geq \gamma$ with asymptotically optimal $\ARL$-$\WADD$ trade-off as $\gamma$ tends to infinity.
	
	\subsection{Unequal Switching Costs}\label{subsect:policy_design_unequal_cost}
	In this subsection, we propose a method to solve the policy design problem in which the switching costs are not all equal. First, we present a proposition regarding the structure of asymptotically optimal solutions of \cref{eqn:QCD_AOSC_formulation}.  
	
	\begin{Proposition}\label{prop:M=1_suffices}
		Suppose $(\tau_{\pi},\pi)$ is an asymptotically optimal solution for \cref{eqn:QCD_AOSC_formulation} with finite window size of at least 1. There exists an open-loop policy $\pi_1$ with window size $W=1$ such that $\AOSC(\pi_1)=\AOSC(\pi)$ and
		\begin{align}
			\WADD(\tau_{\pi_1},\pi_1)&\asymp\WADD(\tau_{\pi},\pi)\quad\text{as $\gamma\to\infty$.}
		\end{align}
	\end{Proposition}
	\begin{IEEEproof}
		See \cref{sec:AppProp5}.
	\end{IEEEproof}
	
	From \cref{prop:M=1_suffices}, the open-loop policy design problem for window size $W\in\bbN$ can be reduced to a problem of window size $\min(W,1)$. Thus, we only need to study the cases where $W=0$ or $W=1$. 
	It should be noted that \cref{prop:M=1_suffices} does not hold when $W=\infty$. An example is provided in \cref{sec:AppExp1}. \cref{prop:M=1_suffices} holds primarily because of the Markovian structure of the sampling policy together with the form of $\AOSC$ function.
	
	In the following, we present algorithms to solve \cref{eqn:optimize_sampling_policy_asymptotic} for each of these cases.
	
	\subsubsection{Window size \texorpdfstring{$W=0$}{W is 0}} 
	Using a similar argument from \cref{subsect:policy_design_equal_cost}, we can see that for any optimal open-loop policy $\pi=(q,T)$ of window size $W$, we have $\overline{q}=q$ and $T$ has only one recurrence class. Thus,  we have 
	\begin{align}
		I_{m,\pi}=\sum_{a\in \calA} q(a) \KLD{p_m(\,\cdot\mid a)}{p_0(\,\cdot\mid a)},\\
		\AOSC(\pi)=\sum_{a\in \calA}\sum_{b\in \calA} q(a)q(b)\calC[a,b],
	\end{align}
	and \cref{eqn:optimize_sampling_policy_asymptotic} becomes
	\begin{equation}\label[Problem]{eqn:policy_design_problem_m=0_not_equal}
	\begin{aligned}
	& \min_{q} & & \max_{1\leq m\leq M}\parens*{\sum_{a\in \mathcal{A}} q(a) \KLD{p_m(\,\cdot\mid a)}{p_0(\,\cdot\mid a)}}^{-1}\\
	& \st
	& & \sum_{a\in\mathcal{A}}q(a)=1, \\
	& & & q(a)\geq 0 \quad \text{for all $a\in\mathcal{A}$},\\
	& & &\sum_{a\in \calA}\sum_{b\in \calA} q(a)q(b)\calC[a,b]\leq \alpha_{\text{AOSC}}.
	\end{aligned}
	\end{equation}
	Using the same argument from \cref{eqn:policy_design_problem_m=0}, we can introduce a new variable $z$ to obtain a linear cost function
	\begin{equation}\label[Problem]{eqn:policy_design_problem_m=0_not_equal_linearize_objective}
	\begin{aligned}
	& \min_{q,z} & &z\\
	& \st
	& & \sum_{a\in\mathcal{A}}q(a)=1,\\
	& & & q(a)>0 \quad \text{for all $a\in\mathcal{A}$},\\
	& & &\sum_{a\in \calA}\sum_{b\in \calA} q(a)q(b)\calC[a,b]\leq \alpha_{\text{AOSC}},\\
	& & &\sum_{a\in \mathcal{A}} q(a) D\left(p_0(\,\cdot\mid a)\ \|\ p_m(\,\cdot\mid a)\right)+z\geq 0\\
	& & &\text{for all $m\in\{1,\ldots,M\}$}.
	\end{aligned}
	\end{equation}
	This is a \gls{QCQP}, and we may assume that $\calC$ is symmetric without loss of generality. However, without additional assumptions, the problem is NP-hard. 
	
	First, we discuss some special cases where the global optimum can be obtained easily. When $\calC$ is positive semi-definite, \cref{eqn:policy_design_problem_m=0_not_equal_linearize_objective} is a convex programming problem. A convex program solver\cite{cvx,gb08} can be used to obtain globally optimal solutions. For the case where there are only cost of observations rather than cost of switching (i.e., $\calC[a,b]=h(b)$ for some function $h:\calA\to \mathbb{R}$), the quadratic constraint in \cref{eqn:policy_design_problem_m=0_not_equal_linearize_objective} reduces to
	\begin{align}\label{eqn:linear_AOSC_constraint}
		\sum_{b\in\calA}q(b)h(b)\leq \alpha_{\text{AOSC}}.
	\end{align} In this case, \cref{eqn:policy_design_problem_m=0_not_equal_linearize_objective} becomes a linear programming problem, which can be solved by a linear program solver\cite{cvx}.
	
	For the general case, we use the IRM algorithm \cite{sun2017rank} to obtain a locally optimal solution. 
	In order to apply the IRM algorithm, we have to first convert \cref{eqn:policy_design_problem_m=0_not_equal_linearize_objective} into a rank-constrained convex optimization problem. We first rewrite the quadratic constraint as
		\begin{align*}
			\sum_{a\in \calA}\sum_{b\in \calA} q(a)q(b)\calC[a,b]=q^T\calC q=\Tr(\calC qq^T)\leq \alpha_{\text{AOSC}}.
		\end{align*}
		By introducing a new  $|\calA|\times |\calA|$ variable $Q$ such that $Q=qq^T$, the quadratic constraint becomes a linear constraint $\Tr(\calC Q)\leq \alpha_{\text{AOSC}}$. To ensure that $Q=qq^T$ holds, we require that $Q\succeq 0$, $Q\mathds{1}=q$ and $\rank(Q)=1$ where $\mathds{1}$ is a $|\calA|\times 1$ vector of ones. Hence, \cref{eqn:policy_design_problem_m=0_not_equal_linearize_objective} is equivalent to 
	\begin{equation}\label[Problem]{eqn:policy_design_problem_m=0_not_equal_SDR}
	\begin{aligned}
	& \min_{Q,q,z} & &z\\
	& \st
	& & \sum_{a\in\mathcal{A}}q(a)=1,\ q(a)>0 \quad \text{for all $a\in\mathcal{A}$},\\
	& & &\Tr(\calC Q)\leq \alpha_{\text{AOSC}},\\
	& & & Q\mathds{1}=q,\  Q\succeq 0\quad \text{and}\quad\rank(Q)=1,\\
	& & &\sum_{a\in \mathcal{A}} q(a) D\left(p_0(\,\cdot\mid a)\ \|\ p_m(\,\cdot\mid a)\right)+z\geq 0\\
	& & &\quad\quad\text{for all $m\in\{1,\ldots,M\}$}.
	\end{aligned}
	\end{equation}
	
	We note that \cref{eqn:policy_design_problem_m=0_not_equal_SDR} becomes a convex programming problem when the constraint $\rank(Q)=1$ is ignored. We are now ready to present the IRM algorithm \cite{sun2017rank}. Fix $\omega>1$.
	
	First, we solve the convex problem 
	\begin{equation}\label[Problem]{eqn:policy_design_problem_m=0_not_equal_SDR_relaxed_t=0}
	\begin{aligned}
	& \min_{Q,q,z} & &z\\
	& \st
	& & \sum_{a\in\mathcal{A}}q(a)=1,\\
	& & & q(a)>0 \quad \text{for all $a\in\mathcal{A}$},\\
	& & &\Tr(\calC Q)\leq \alpha_{\text{AOSC}},\ Q\mathds{1}=q,\ Q\succeq 0,\\
	& & &\sum_{a\in \mathcal{A}} q(a) D\left(p_0(\,\cdot\mid a)\ \|\ p_m(\,\cdot\mid a)\right)+z\geq 0\\
	& & &\quad\quad\text{for all $m\in\{1,\ldots,M\}$},
	\end{aligned}
	\end{equation}
	to obtain a solution $(Q_0,q_0,z_0)$ and let $Q_0=VDV^T$ be the eigen-decomposition of $Q_0$. Let $V_0$ be the eigenvectors corresponding to the $|\calA|-1$ smallest eigenvalues of $Q_0$. 
	
	At each step $k$, we solve the following convex problem:
	\begin{equation}\label[Problem]{eqn:policy_design_problem_m=0_not_equal_SDR_relaxed_t=k}
	\begin{aligned}
	& \min_{Q,q,z,r} & &z+\omega^kr\\
	& \st
	& & \sum_{a\in\mathcal{A}}q(a)=1,\\
	& & & q(a)>0 \quad \text{for all $a\in\mathcal{A}$},\\
	& & &\Tr(\calC Q)\leq \alpha_{\text{AOSC}},\ Q\mathds{1}=q,\ Q\succeq 0,\\
	& & & rI_{n-1}-V_{k-1}^TQV_{k-1}\succeq 0,\\
	& & &\sum_{a\in \mathcal{A}} q(a) D\left(p_0(\,\cdot\mid a)\ \|\ p_m(\,\cdot\mid a)\right)+z\geq 0\\
	& & &\quad\quad\text{for all $m\in\{1,\ldots,M\}$},
	\end{aligned}
	\end{equation}
	to obtain a solution $(Q_k,q_k,z_k,r_k)$ and let $V_k$ be the eigenvectors corresponding to the $|\calA|-1$ smallest eigenvalues of $Q_k$. We iterate until $r_k<\epsilon$ , where $\epsilon$ is a small threshold chosen as a stopping criterion. Following similar methods from \cite{sun2017rank}, it can be shown that $r_k\to 0$ at a linear rate and that $q_k$ converges to a locally optimal solution for \cref{eqn:policy_design_problem_m=0_not_equal_SDR} if \cref{eqn:policy_design_problem_m=0_not_equal_SDR_relaxed_t=k} is feasible for all $k$.
	
	\subsubsection{Window size \texorpdfstring{$W=1$}{W is 1}}
	Unlike the case when $W=0$, not every distribution $q$ is a stationary distribution of $T$. Furthermore, when $W>0$, it is possible that more than one recurrence class exists. In this case, \cref{eqn:optimize_sampling_policy_asymptotic} becomes
	\begin{equation}\label[Problem]{eqn:policy_design_problem_m=1}
		\begin{aligned}
			& \min_{T,q} & & \max_{1\leq m\leq M} (\sum_{a\in \mathcal{A}} q(a) D\left(p_0(\,\cdot\mid a)\ \|\ p_m(\,\cdot\mid a)\right))^{-1}\\
			& \st
			& & \sum_{a\in\mathcal{A}}q(a)=1, \ q^TT=q^T,\\
			& & & q(a)>0 \quad \text{for all $a\in\mathcal{A}$},\\
			& & &\sum_{b\in\mathcal{A}}T[a,b]=1 \quad \text{for all $a\in\mathcal{A}$},\\
			& & &0\leq T[a,b]\leq 1 \quad \text{for all $a,b\in\mathcal{A}$},\\
			& & &\sum_{b\in\mathcal{A}}\sum_{a\in\calA}T[a,b]\calC[a,b]q(a)\leq \alpha_{\text{AOSC}},\\
			& &&\supp(q) \subseteq \text{ one recurrence class}.
		\end{aligned}
	\end{equation}
	Using the same argument from \cref{eqn:policy_design_problem_m=0}, \cref{eqn:policy_design_problem_m=1} is equivalent to 
	\begin{equation}\label[Problem]{eqn:qcqp_policy_design_problem_m=1}
	\begin{aligned}
	& \min_{T,q,z} & &z\\
	& \st
	& & \sum_{a\in\mathcal{A}}q(a)=1, \\
	& & & q(a)>0 \quad \text{for all $a\in\mathcal{A}$}\\
	& & &\sum_{b\in\mathcal{A}}T[a,b]=1 \quad \text{for all $a\in\mathcal{A}$}\\
	& & &0\leq T[a,b]\leq 1 \quad \text{for all $a,b\in\mathcal{A}$}\\
	& & &q^TT=q^T,\\
	& & &\sum_{b\in\mathcal{A}}\sum_{a\in\calA}T[a,b]\calC[a,b]q(a)\leq \alpha_{\text{AOSC}},\\
	& &&\supp(q) \subseteq \text{ one recurrence class},\\
	& & &\sum_{a\in \mathcal{A}} q(a) D\left(p_m(\,\cdot\mid a)\ \|\ p_0(\,\cdot\mid a)\right)+z\geq 0,\\
	& & &\quad\quad\text{for all $m\in\{1,\ldots,M\}$.}
	\end{aligned}
	\end{equation}
	\cref{eqn:qcqp_policy_design_problem_m=1} has two quadratic constraints $Tq=q$ and $\sum_{b\in\mathcal{A}}\sum_{a\in\calA}T[a,b]\calC[a,b]q(b)\leq \alpha_{\text{AOSC}}$. Thus, it is a \gls{QCQP} with an additional combinatorial constraint that  $\supp(q)$ is contained in a single recurrence class of $T$. Even without the combinatorial constraint, finding the global optimal for \cref{eqn:qcqp_policy_design_problem_m=1} would be difficult without additional assumptions on $\calC$.
	
	By considering a similar construction used in the proof of \cref{prop:M=1_suffices}, we show that any open-loop policy $\pi_1=(q,T)$ with window size $1$ can be expressed as an open-loop policy $\pi_2=(q_2,T_2)$ of window size $2$ satisfying the following:
\begin{align}
			\sum_{b\in \calA} q_2(a,b)=\sum_{c\in \calA} q_2(c,a)\quad\text{for any $a\in\calA$},\label{eqn:pi1topi2_1}\\
			T_2[(a,b),(c,d)]=\begin{cases}
				\frac{q_2(c,d)}{\sum_{d\in\calA}q_2(c,d)}\quad\text{if $b=c$,}\\
				0\quad\text{otherwise.}
			\end{cases}\label{eqn:pi1topi2_2}
	\end{align}
	Furthermore, any open-loop policy $\pi_2=(q_2,T_2)$ with window size $2$ that satisfies \cref{eqn:pi1topi2_1,eqn:pi1topi2_2} can be expressed as an open-loop policy $\pi_1$ of window size $1$. When we consider an open-loop policy $\pi_2=(q_2,T_2)$ with window size $2$ that satisfies \cref{eqn:pi1topi2_1,eqn:pi1topi2_2}, the constraint $T_2q_2=q_2$ becomes automatically satisfied and the $\AOSC$ constraint is linearized to $\sum_{(a,b)\in\mathcal{A}^2}\calC[a,b]q(a,b)\leq \alpha_{\text{AOSC}}$. Hence, \cref{eqn:qcqp_policy_design_problem_m=1} is equivalent to the following problem
	\begin{equation}\label[Problem]{eqn:lp_policy_design_problem_m=1}
		\begin{aligned}
		& \min_{q,z}& &z\\
		& \st
		& & \sum_{(a,b)\in\mathcal{A}^2}q(a,b)=1, \\ 
		& & & q(a,b)>0\quad \text{for all $a,b\in\mathcal{A}$},\\
		& & &\sum_{(a,b)\in\mathcal{A}^2}\calC[a,b]q(a,b)\leq \alpha_{\text{AOSC}},\\
		& & &\sum_{b\in \calA} q(a,b)=\sum_{c\in \calA} q(c,a)\quad\text{for all $a\in\mathcal{A}$}, \\
		& & &\supp(q) \subseteq \text{ one recurrence class},\\
		& & &\sum_{(a,b)\in A^2} q(a,b) D\left(p_m(\,\cdot\mid b)\ \|\ p_0(\,\cdot\mid b)\right)\geq -z,\\
		& & &\quad\quad\text{for all $m\in\{1,\ldots,M\}$, $a,b\in\mathcal{A}$},
		\end{aligned}\end{equation}
	which is a linear programming problem with an additional combinatorial constraint. In order to handle the combinatorial constraint, we can solve \cref{eqn:lp_policy_design_problem_m=1} without the combinatorial constraint $\supp(q) \subseteq \text{ one recurrence class}$. The solution obtained by solving \cref{eqn:lp_policy_design_problem_m=1} without the combinatorial constraint will not have support in one recurrence class in general. However, if the solution has support in one recurrence class, it is an optimal solution to \cref{eqn:lp_policy_design_problem_m=1}. Alternatively, we can select a sufficiently small $\epsilon>0$, and require that $q(a,b)>\epsilon$ for all $a,b\in\mathcal{A}$. The new constraint ensures that there is only one recurrence class for any feasible open-loop policy and thus, the recurrence class constraint is satisfied. The relaxed problem becomes
	\begin{equation}\label[Problem]{eqn:lp_policy_design_problem_m=1_relaxed}
	\begin{aligned}
	&\min_{q,z} & & z\\
	& \st
	& & \sum_{(a,b)\in\mathcal{A}^2}q(a,b)=1,\\
	& & & q(a,b)>\epsilon\ \quad\text{for all $a,b\in\mathcal{A}$},\\
	& & & \sum_{(a,b)\in\mathcal{A}^2}\calC[a,b]q(a,b)\leq \alpha_{\AOSC},\\
	& & & \sum_{a\in \calA} q(a,b)=\sum_{c\in \calA} q(c,b)\ \text{for $b\in\calA$,} \\
	& & & \sum_{(a,b)\in A^2} q(a,b) D\left(p_m(\,\cdot\mid b)\ \|\ p_0(\,\cdot\mid b)\right)\geq -z,\\
	& & & \quad\quad \text{for all $m\in\{1,\ldots,M\}$},
	\end{aligned}
	\end{equation}
	which can be solved by a linear program solver\cite{cvx}.
\section{Observation-dependent Sampling Policy}\label{sec:closed_loop_policy_design}
In the discussion of the preceding sections, we focused on open-loop policies for which the actions used to obtain the observations are independent of the observations obtained during test time. For some applications, the physical infrastructure allows real-time control of the actions based on the observations obtained during test time. In the case when the post-change distribution is known, this control can used to select actions that minimizes the $\WADD$ for that particular post-change distribution. Thus, the additional control can be used to further lower the $\WADD$ as compared to the open-loop policy. In this section, we propose an observation-dependent sampling policy for \cref{eqn:QCD_AOSC_formulation} that uses this additional control to improve the performance obtained using the open-loop policies. 
\subsection{Policy Design}
The design of the observation-dependent sampling policy is motivated by different considerations in the pre- and post-change regimes. We consider a two-stage sampling policy where sampling policies used in the first and second stage is designed to satisfy the considerations in the pre- and post-change regime, respectively. 

Since the observations obtained in the pre-change regime are i.i.d. and contains no information about post-change distribution $p_m$, observation-dependent policies are unable to make use of these observations to reduce the ASOC or improve the $\ARL$-$\WADD$ trade-off rate. Hence, we only need to consider open-loop policies for the first stage. In the pre-change regime, the main consideration in designing the sampling policy is to satisfy the $\AOSC$ constraint while maximizing distinguishability of the pre- and post-change distributions. The open-loop sampling policy that solves \cref{eqn:optimize_sampling_policy_asymptotic} is an open-loop sampling policy that satisfies these considerations and hence, is chosen to be the sampling policy for the first stage of the two-stage sampling policy.

On the other hand, the main consideration for the design of the sampling policy in the post-change regime is to optimize the $\ARL$-$\WADD$ trade-off rate. The observations obtained in the post-change regime allow us to make meaningful inference of the post-change distribution $p_m$ and hence are important for selecting actions that optimize the $\ARL$-$\WADD$ trade-off rate. One way to approach this coupled decision problem of inferring $p_m$ and selecting the action is the Chernoff rule \cite{chernoff1959sequential}. For our QCD problem, applying the Chernoff rule selects the current action that maximizes the KL divergence  $\KLD{p_{m_{\textbf{ML}}}(\,\cdot\mid A)}{p_0(\,\cdot\mid A)}$ where $m_{\textbf{ML}}$ is the maximum likelihood (ML) decision about the post-change distribution. To ensure that the ML decision $m_{\textbf{ML}}$ converges to the true $m$ in finite time, the Chernoff rule\cite{chernoff1959sequential} requires the additional assumption that for all actions $A$ and for each $m\in\{1,\ldots,M\}$ , we have $\KLD{p_m(\,\cdot\mid A)}{p_0(\,\cdot\mid A)}>0$. In \cite{nitinawarat2013controlled}, it is shown that this assumption is not required if a randomized selection rule is performed on the actions $\calA$ at exponentially spaced time instances for the application of sequential multi-hypothesis testing. Thus, for the second-stage of the observation-dependent sampling policy, we select the current action to maximize $\KLD{p_{m_{\textbf{ML}}}(\,\cdot\mid A)}{p_0(\,\cdot\mid A)}$ and randomly select an action at exponentially spaced intervals. 

We propose a $2$-threshold stopping time $\tau_{\widehat{\pi}}$ that uses a $2$-stage sampling policy $\widehat{\pi}$ to select the actions while using the GLR CuSum.  We let $A^*_m$ be the action that maximizes the KL divergence $\KLD{p_m(\,\cdot\mid A)}{p_0(\,\cdot\mid A)}$, i.e.,
\begin{align*}
A^*_m=\argmax_{A\in\calA} \KLD{p_m(\,\cdot\mid A)}{p_0(\,\cdot\mid A)}
\end{align*}
for $m=1,\ldots,M$, $\pi^\dagger=(q,T)$ be an open-loop policy of window size $2$ that solves \cref{eqn:optimize_sampling_policy_asymptotic}, $\kappa,\gamma$ be thresholds such that $0<\kappa<\gamma$ and base $b>1$ for determining exponentially spaced intervals. The proposed sampling policy $\widehat{\pi}$ is implemented as follows. We run the GLR CuSum stopping time. We use the $2$-stage sampling policy to determine the actions used to obtain the observations. We initialize by using actions generated by the open-loop policy $\pi^\dagger=(q,T)$ in the first stage. Sampling using $\pi^\dagger$ continues until time $t_{\kappa}$ which is the first time the GLR CuSum statistic $S(\widehat{\pi},t)$ exceeds the first threshold $\kappa$.  We transit into the second stage and now decide the actions by selecting an action $A^*_{m_{\text{ML}}}$ where $m_{\text{ML}}$ is the ML decision of the post-change distribution and randomly selecting from $\calA$ at time instances $t=t_{\kappa}+\ceil{b^l}$ where $l\in\mathbb{N}$. We continue sampling until either the GLR CuSum $S(\widehat{\pi},t)$ falls below the first threshold $\kappa$, and we return to sampling using the open-loop policy $\pi^\dagger$, or $S(\widehat{\pi},t)$ exceeds the second threshold $\gamma$ and a change is declared. Our proposed algorithm is described in detail in \cref{alg:closedloop}.  
\begin{algorithm}[!htbp]
	\caption{$2$-threshold stopping time $\tau_{\widehat{\pi}}$ given an open-loop policy $\pi^\dagger=(q,T)$ of window size $2$, thresholds $\kappa,\gamma$ with $0<\kappa<\gamma$ and $b>1$.}\label{alg:closedloop}
	\begin{algorithmic}
		\State $\textbf{Initialize:}$
		\State Set $t=2$, $l=0$, $t_{\kappa}=0$
		\State Sample $\calA^2$ using $q$ to obtain the first two actions $a_1,a_2$.
		\State Compute $S(\widehat{\pi},t)$ according to \cref{eqn:parallel_CuSum_stoptime}.
		\While {$S(\widehat{\pi},t)\leq\gamma$}
		\If {$S(\widehat{\pi},t)\leq \kappa$}
		\If {$S(\widehat{\pi},t-1)\leq \kappa$}
		\State Obtain action $a_{t+1}$ according $\pi^\dagger$.
		\Else
		\State Sample $\calA^2$ using $q$ to obtain the next two actions\\\quad\quad\quad\quad\quad $a_{t+1},a_{t+2}$ for policy $\pi^\dagger$.
		\EndIf
		\Else
		\If{$S(\widehat{\pi},t-1)\leq \kappa$}
		\State $t_\kappa:=t$, 	$l:=0$.
		\EndIf
		\If{$t=t_{\kappa}+\ceil{b^l}$}
		\State Select $a_{t+1}$ from $\calA$ randomly with uniform\\\quad\quad\quad\quad\ distribution.
		\State $l:=l+1$.
		\Else
		\State $a_{t+1}:=A^*_m$ where $\displaystyle m=\argmax_{1\leq m\leq M}S_m(\widehat{\pi},t-1)$.
		\EndIf
		\EndIf
		\State $t:= t+1$
		\State Update $S(\widehat{\pi},t)$ according to \cref{eqn:parallel_CuSum_update}.
		\EndWhile
		\State Stop sampling and declare a change $\tau_{\widehat{\pi}}=t$.
	\end{algorithmic}
\end{algorithm}

Next, we show that by properly selecting a threshold $\kappa$, we are able to control the $\AOSC$ of the policy $\widehat{\pi}$.
\begin{Proposition}\label{prop:ASOC}
	Let $\alpha_{\AOSC}>0$ be such that \cref{eqn:optimize_sampling_policy_asymptotic} is feasible with a strictly positive optimal value and optimal solution $\pi^\dagger$. For any $\theta>0$, there exists $\kappa>0$ such that for all $\gamma>\kappa$ and $b>1$, the $2$-threshold stopping time $\tau_{\widehat{\pi}}$ satisfies
	\begin{align}
	\AOSC(\widehat{\pi})\leq \alpha_{\AOSC}+\theta.
	\end{align}
\end{Proposition}
\begin{IEEEproof}
Let $\psi_1$ be the proportion of action-switches that is generated by $\pi^\dagger$ and $\psi_2$ be the proportion of action-switches that is not generated by $\pi^\dagger$ for any sequence of actions $\{a_t\}_{t\in\mathbb{N}}$ generated by the casual sampling policy $\widehat{\pi}$. For any sequence of actions $\{a_t\}_{t\in\mathbb{N}}$, the sum of log-likelihood ratios $\sum_{i=1}^t\log\frac{p_m(\,X_t \mid a_t)}{p_0(\,X_t \mid a_t)}$ has a negative drift under $\P_\infty$. Thus, the proportion $\psi_2$ approaches zero as $\kappa\to\infty$. We select a sufficiently large threshold $\kappa$ such that $\psi_2 \max_{i,j}\mathcal{C}[i,j]<\theta$. Since $\AOSC(\pi)\leq \psi_1\alpha_{\text{AOSC}}+\psi_2 \max_{i,j}\calC[i,j]$, we have $\AOSC(\pi)\leq \alpha_{\AOSC}+\theta$ and the proof is complete.
\end{IEEEproof}

Using a similar argument to \cref{prop:Information_number_WADD}, it can be shown that $\ARL(\tau_{\widehat{\pi}},\widehat{\pi})\geq \gamma$.
\begin{Proposition}\label{prop:2thresARL}
	Let $\gamma>0$. For all $0<\kappa<\gamma$,  stopping time $\tau_{\widehat{\pi}}$ satisfies
	\begin{align}
	\ARL(\tau_{\widehat{\pi}},\widehat{\pi})\geq \gamma.
	\end{align}
\end{Proposition} However, we are unable to obtain the upper-bound for the $\WADD(\tau_{\widehat{\pi}},\widehat{\pi})$ due to theoretical difficulties in deriving similar results for \eqref{eqn:information_number_as_convergence}, \eqref{eqn:information_number_convergence_in_prob} and \eqref{eqn:information_number_convergence_in_prob_lower_bound}. Our simulations indicate that the $\ARL$-$\WADD$ tradeoff rate approaches the trade-off rate achieved by the GLR CuSum when the post-change distribution is known.
\subsection{Choice of parameters $\gamma,\kappa$ and $b$}
By \cref{prop:2thresARL}, the $\ARL$ of the $2$-threshold stopping time is at least $\gamma$. Hence, during implementation, we set $\gamma$ to be the $\ARL$ requirements required by the application. From the proof of \cref{prop:ASOC}, the $\AOSC$ is a decreasing function of $\kappa$. The parameter $\kappa$ is can be selected to satisfy $\AOSC(\widehat{\pi})\leq \alpha_{\AOSC}+\theta$ using 
Monte-Carlo approximations of $\AOSC(\widehat{\pi})$ and a grid search between $0$ and $\gamma$. As $\gamma$ is set by the average runtime to false alarm of the application, it is usually large and a suitable $\kappa$ can be found within $0$ and $\gamma$. If a suitable $\kappa$ cannot be found, $\gamma$ can be increased so that $0<\kappa<\gamma$ at the expense of a larger $\WADD$. Finally, $b$ controls the exponentially spaced time instances, after the test statistic $S(\widehat{\pi},t-1)$ exceeds $\kappa$, at which the actions are randomly selected. The purpose of randomly selecting the action is to prevent the ML decision about the post-change distribution from being stuck at a wrong decision. Using a large $b$ would result in a large $\WADD$ in the event that the ML decision becomes stuck at a wrong decision. On the other hand using a small $b$ would result in frequent deviation from the action that optimizes the $\ARL$-$\WADD$ trade-off rate and hence a degradation of performance when the $\gamma$ is small. For the numerical simulations performed in this paper, we set $b=2$.

	\section{Numerical Results}\label{sec:results}
	In this section, we present numerical results to illustrate the relationship between the asymptotic ARL-WADD trade-off rate and the AOSC of the open-loop policy. We also compare our proposed sampling policy, which takes observation-switching costs constraints into account, against observation-dependent polices with only sampling cost constraints.
	\subsection{Applications to graph signals}
\begin{figure}[H]
		\centering
		\centerline{\includegraphics[width=0.9\columnwidth]{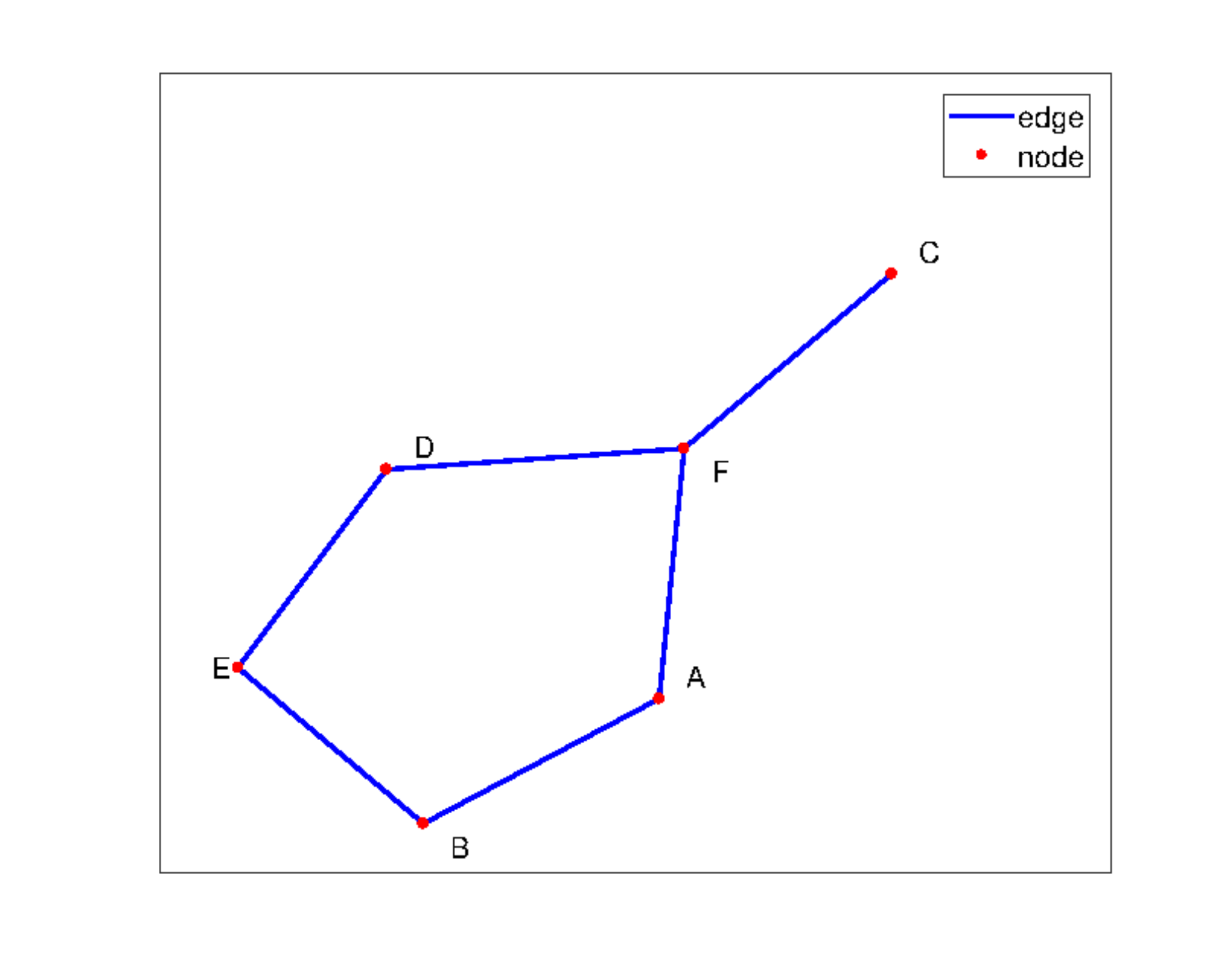}}
		\caption{Graph $\calG$ generated using the Erd\H{o}s-R\'enyi (ER) random graph model with $N=6$ nodes and probability of an edge $p=0.2$.}
		\label{fig:graph}
	\end{figure}

	In this subsection, we consider the QCD problem with $\AOSC$ based on partially observed graph signals under the various conditions discussed in this paper. We consider the problem of quickest detection of a rogue node in a graph. We assume that our graph $\mathcal{G}$ is a connected graph with $N$ nodes. 
	We model the graph signal\cite{dong2016learning} in the pre-change regime with a zero-mean Gaussian distribution with covariance $\Sigma=L^\dagger+\eta I$, where $L^\dagger$ is the psuedo-inverse of the graph Laplacian matrix $L$, $I$ is an $N\times N$ identity matrix and $\eta^2$ is the noise power in the graph signal. For all the simulations in this section, the noise power $\eta^2$ is set at $0.01$.
	 Thus, in the pre-change regime, we have
	\begin{align}
		X_t\sim p_0&=\mathcal{N}(\mathbf{0},\Sigma).
	\end{align}
	For the post-change regime, we assume that the signal obtained at the rogue node follows the same distribution as the pre-change distribution. However, this signal becomes independent of signals obtained from the rest of the graph. Thus, in the post-change regime, we have
	\begin{align}
		X_t\sim p_m&=\mathcal{N}(\mathbf{0},\Sigma_m)
	\end{align} 
	for $m\in\{1,\ldots,M\}$, where the covariance matrix $\Sigma_m$ is given as
	\begin{align}
		\Sigma_m[i,j]=
		\begin{cases}
			\Sigma[i,j] &\text{for $i\neq m$ and $j\neq m$,}\\
			\Sigma[m,m] &\text{for $i=j= m$,}\\
			0&\text{otherwise.}
		\end{cases}\label{eqn:post_change}
	\end{align}
	The set of actions $\mathcal{A}$ is the set of partial observations where 
	\begin{align}
		&\mathcal{A}=\bigcup_{n=2}^{N-2}\braces*{\mathbf{M}\ |\  \mathbf{M}=[e_{i_1},e_{i_2},...,e_{i_n}]^T,i_1<i_2<\ldots<i_n}.
	\end{align}
	In our experiments, we consider the graph $\calG$, as shown in \figref{fig:graph}. In this case, there is a total of $50$ actions that can be used to observe the signal on the graph $\calG$.

	We consider three possible cases of observation-switching costs $\calC$: 
	\begin{enumerate}
		\item no costs are involved with $\calC^1=\mathbf{0}$;
		\item there is only activation cost where $\calC^2[i,j]$ is the number of nodes observed using action $i$; 
		\item there are both switching and activation costs with $\calC^3=\calC^2+\calC'$ where $\calC'[i,j]$ is the number of elements in the symmetric difference between the set of nodes observed by action $i$ and set of nodes observed by action $j$.
	\end{enumerate} 
	Note that each $\calC^1,\calC^2,\calC^3\in\mathbb{R}^{50\times50}$.
	
	When $\calC=\calC^1$, using results from \cref{subsect:policy_design_equal_cost}, we only need to design an open-loop policy of window size $0$ if we are considering only open-loop policies. By considering the problem for $W=0$, we obtain the optimal open-loop policy $\pi_{1,W=0}$ by solving \cref{eqn:policy_design_problem_m=0_LP}. Similarly, for the cases where $\calC=\calC^2,\calC^3$, using results from \cref{subsect:policy_design_unequal_cost}, we only need to consider the open-loop policies with $W=0,1$ if we are only considering open-loop policies. When $\calC=\calC^2$ and $W=0$, we obtain the optimal open-loop policy $\pi_{2,W=0}$ by solving \cref{eqn:policy_design_problem_m=0_not_equal_linearize_objective} with a linear $\AOSC$ constraint \cref{eqn:linear_AOSC_constraint}.  For the case $\calC=\calC^3$ with $W=0$, we obtain a locally optimal open-loop policy $\pi_{3,W=0}$ by performing the IRM algorithm described in \cref{subsect:policy_design_unequal_cost}.  In the case where $\calC=\calC^2$ or $\calC^3$ with $W=1$, we obtain an approximately optimal open-loop policy $\pi_{2,W=1,\epsilon}$ or $\pi_{3,W=1,\epsilon}$ respectively by solving \cref{eqn:lp_policy_design_problem_m=1_relaxed} with $\epsilon\in\{10^{-4},10^{-5},10^{-6}\}$.
	
	When no costs are involved (i.e., $\calC=\mathbf{0}$), the optimal asymptotic ARL-WADD trade-off rate (cf.\ \cref{def:trade-off_rate}) achieved by the stopping rule and policy $(\tau_{\pi_{1,W=0}},\pi_{1,W=0})$ is $\bar{I}=1.5586$. This is an upper bound for the asymptotic ARL-WADD trade-off rates for all other choices of $\calC$. 
	
\begin{figure}[!htb]
		\centering
		\centerline{\includegraphics[width=\columnwidth]{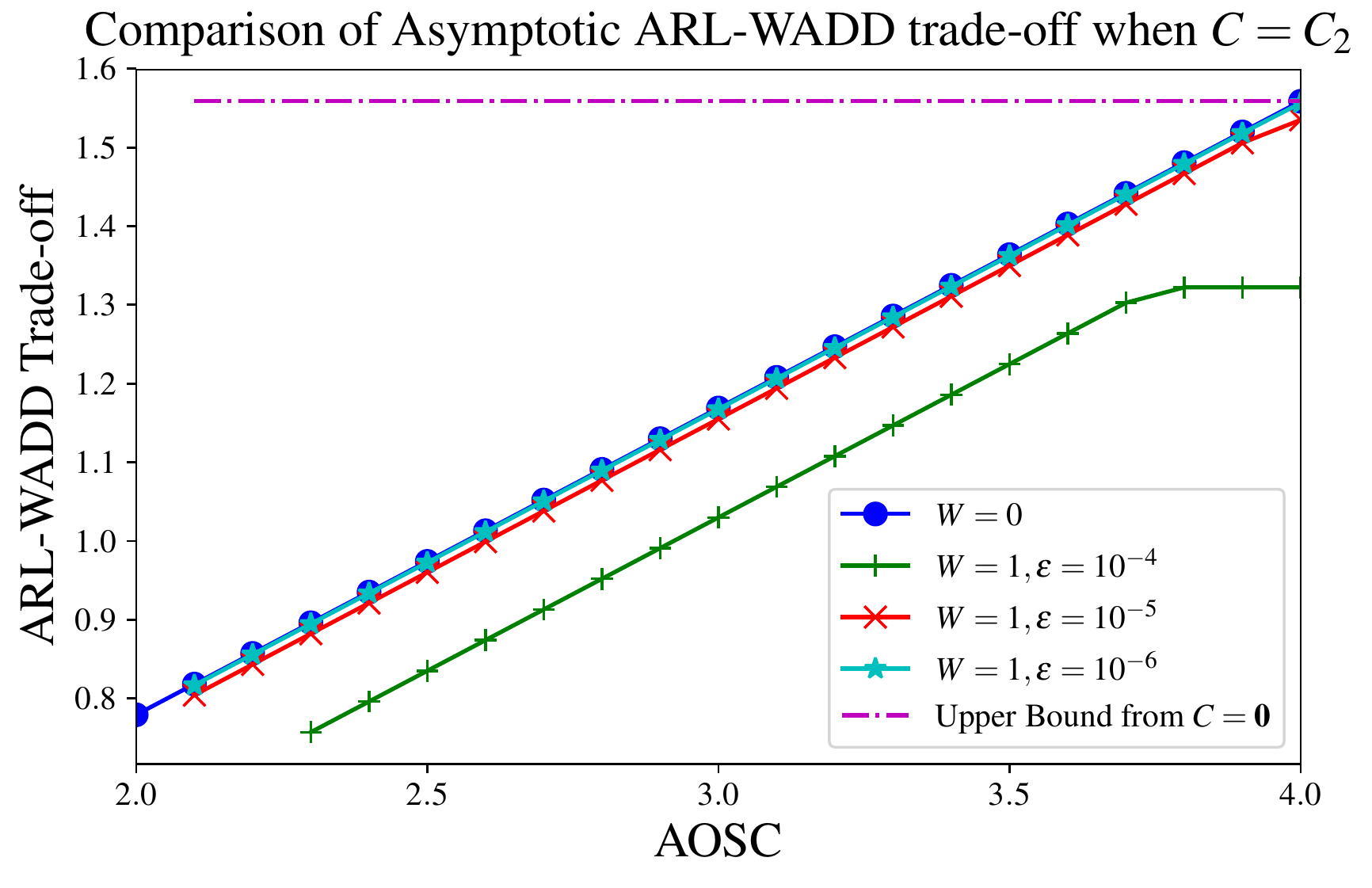}}
		\caption{Comparison of the ARL-WADD trade-off rate for different stopping times when $\calC=\calC^2$. }
		\label{fig:ARLWADDc2}
	\end{figure}

	In \cref{fig:ARLWADDc2}, we compare the asymptotic ARL-WADD trade-off rate of $\pi_{2,W=0}$ and $\pi_{2,W=1,\epsilon}$. First, we observe that the asymptotic $\ARL$-$\WADD$ trade-off is below the upper bound $\bar{I}$ across the range of achievable $\AOSC$. In the case where $\calC=\calC^2$, the $\AOSC$ constraint in \cref{eqn:policy_design_problem_m=0_not_equal_linearize_objective} reduces to
	\begin{align}
		\AOSC(\pi)&=\sum_{b\in \calA}\left(\sum_{a \in\calA}q(a,b)\right)h(b)
	\end{align} where $h(b)$ is the number of nodes observed using action $b$. Using similar arguments from the proof of \cref{prop:M=0_suffices}, we can see that there is an asymptotically optimal trade-off for a policy window of size $W=0$ is equal to the asymptotically optimal trade-off for a policy of window size $W=1$. Hence, the performance between $\pi_{2,W=0}$ and $\pi_{2,W=1,\epsilon}$ becomes more similar as $\epsilon$ tends to zero.
	
\begin{figure}[H]
		\centering
		\centerline{\includegraphics[width=\columnwidth]{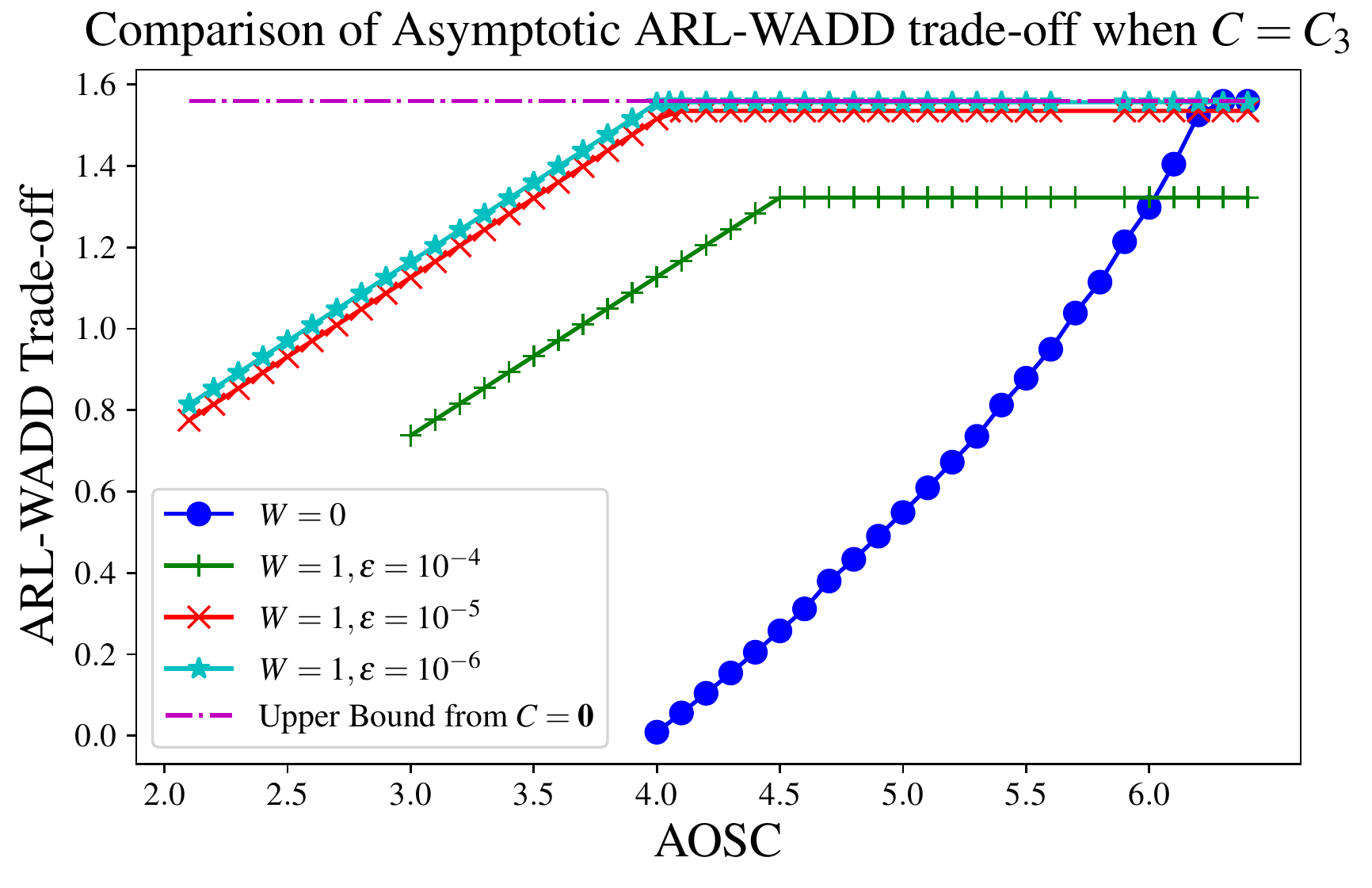}}
		\caption{Comparison of the ARL-WADD trade-off rate for different stopping times when $\calC=\calC^3$. }
		\label{fig:ARLWADDc3}
	\end{figure}

\begin{figure}[!htb]
		\centering
		\centerline{\includegraphics[width=\columnwidth]{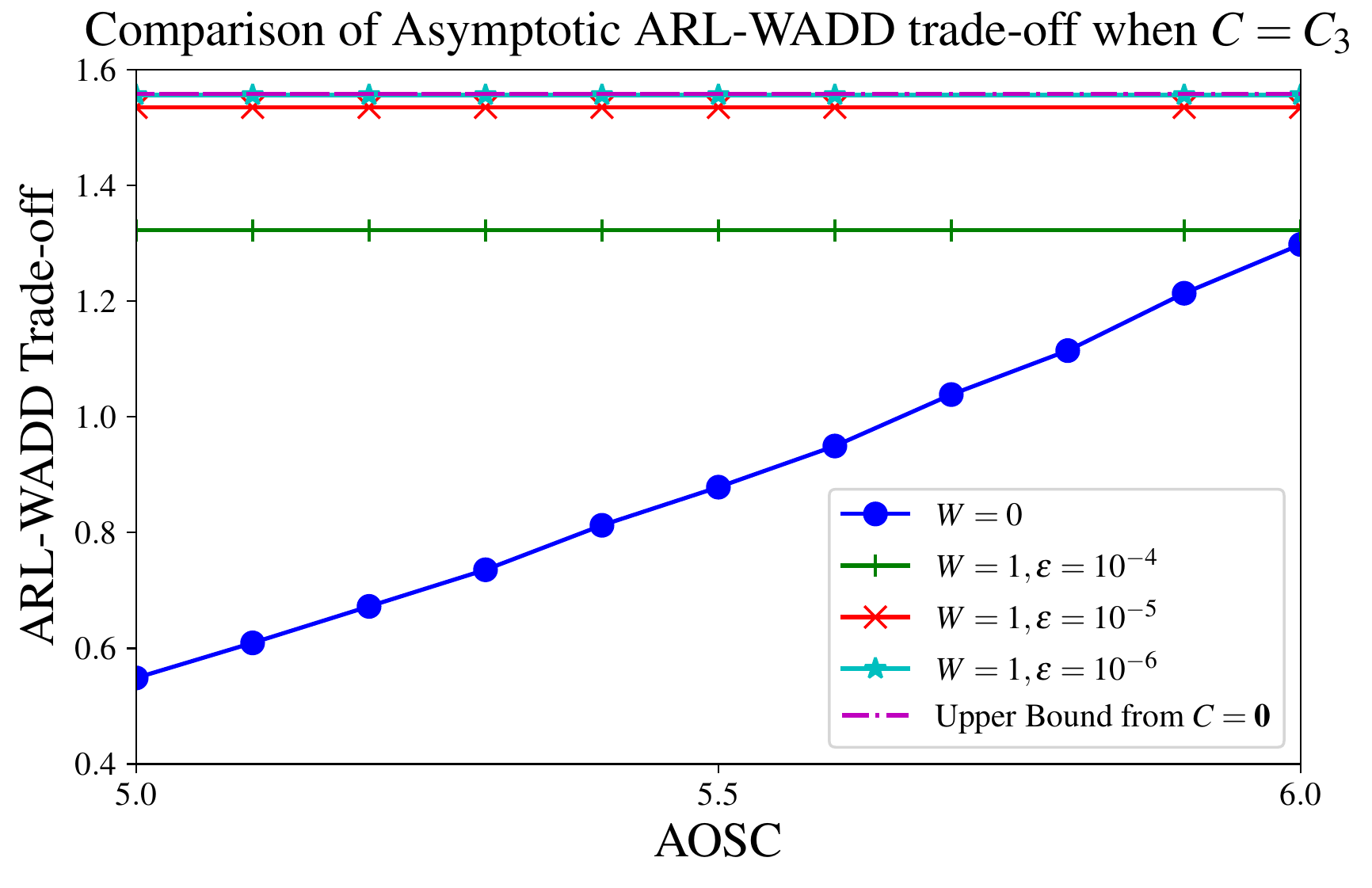}}
		\caption{Zoomed version of \figref{fig:ARLWADDc3} into $\AOSC\geq 5$. }
		\label{fig:ARLWADDc3_zoom}
	\end{figure}

	In \cref{fig:ARLWADDc3,fig:ARLWADDc3_zoom}, we compare the asymptotic ARL-WADD trade-off rate of $\pi_{3,W=0}$ and $\pi_{3,W=1,\epsilon}$. Similarly, we observe that the asymptotic $\ARL$-$\WADD$ trade-off is below the upper bound $\bar{I}$ across the range of achievable AOSC. In this case, we observe that the optimal open-loop policies of window size $W=1$ significantly outperforms the optimal open-loop policies of window size $W=0$. When we are using a policy of window size $W=0$, the lowest $\AOSC$ achievable while \cref{eqn:QCD_AOSC_formulation} remains feasible is about $4$. However, using a policy of window size $W=1$, we are able to reduce $\AOSC$ to about $2$ while \cref{eqn:QCD_AOSC_formulation} remains feasible. These can be seen be comparing the lowest $\AOSC$ achieved by the respective curves in \cref{fig:ARLWADDc3}.
	
		\subsection{Comparison with observation-dependent policies having only sampling cost constraints}\label{subsec:results}
		While there are no existing work on QCD with action switching-costs, for some special cases, existing work in the literature may be relevant. In this subsection, we compare our proposed GLR CuSum with the GDECuSum proposed in \cite{banerjee15} and discuss their strengths and weaknesses. The GDECuSum is proposed in which an on-off observation control is used to take into account the cost of observations for the purpose of QCD. The main difference between the GDECuSum and our proposed method is that our proposed method uses action switching costs while the GDECuSum only takes observation costs into account. In the next two simulations, we demonstrate the differences between the two methods by using different switching cost matrix settings. We use the same distributions and parameters for the GDECuSum as \cite{banerjee15}. The signal is generated by a pre-change distribution with pdfs $p_0=\calN(0,1)$ and $4$ possible post-change distributions with pdf $p_1=\calN(0.4,1)$, $p_2=\calN(0.6,1)$, $p_3=\calN(0.8,1)$ and $p_4=\calN(1,1)$. The parameters for the GDECuSum are $\mu=0.08$ and $h=\infty$. Using these parameters, the GDECuSum achieves $\ARL-\WADD$ trade-off of $0.08$ and a Pre-change Duty Cycle (PDC) of $0.5$ which means that only $50\%$ of the samples are observed under the pre-change regime. At any time instance $t$, we let $A_t=1$ denote the action where the $t$-th sample is observed and $A_t=2$ denote the action where the $t$-th sample is skipped.
		
\begin{figure}[!htb]
			\centering
			\centerline{\includegraphics[width=\columnwidth]{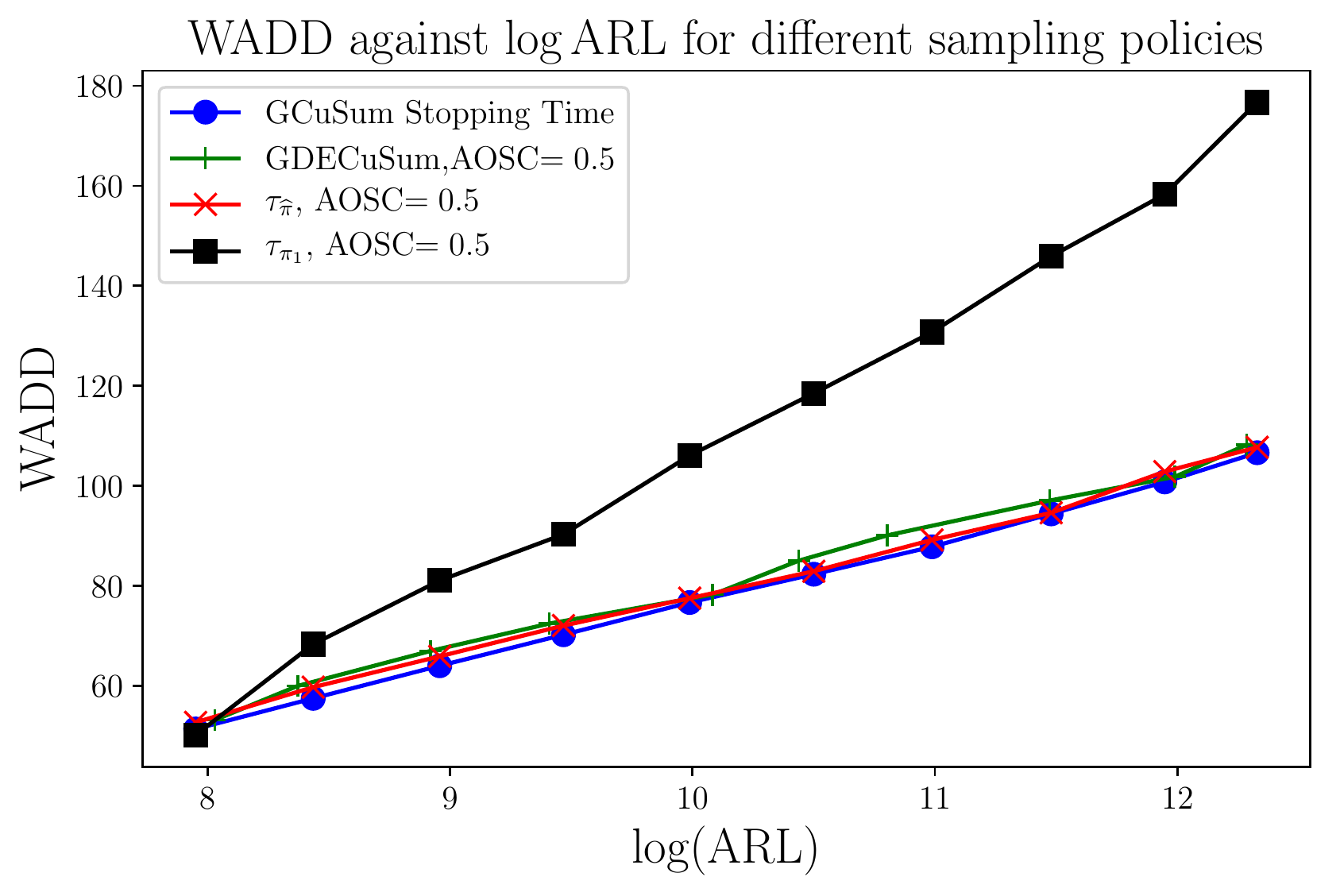}}
			\caption{Comparison of the ARL-WADD trade-off rate for the GDECuSum and the GLR CuSum with the 2-threshold stopping time when $\calC=\calC^4$. }
			\label{fig:ARLWADDc4}
		\end{figure}

		In the first set of simulations, we use the observation-switching cost matrix
	$
			\calC^4=\begin{pmatrix}
				1 & 0\\
				1&	0
			\end{pmatrix}.
		$
		In this case, the $\AOSC$ is equal to the PDC as the $\AOSC$ measures the average proportion of samples observed in the pre-change regime. We solve \cref{eqn:lp_policy_design_problem_m=1_relaxed} with $\AOSC=\text{PDC}=0.5$ and $\epsilon=10^{-6}$ to obtain the $\ARL-\WADD$ trade-off of $0.04$ for our proposed GLR CuSum stopping time with an open-loop policy of window size $1$. Thus, the GDECuSum which uses an observation-dependent policy out-performs our proposed GLR CuSum as expected when we use an open-loop policy of window size $1$. In \figref{fig:ARLWADDc4}, we compare the generalised CuSum (GCuSum) stopping time, GDECuSum\cite{banerjee15} with our proposed GLR CuSum with the 2-threshold stopping time $\tau_{\widehat{\pi}}$ with observation-dependent sampling policy $\widehat{\pi}$. The 2-threshold stopping time is determined by letting $\pi^\dagger$ to be the optimal solution of \cref{eqn:lp_policy_design_problem_m=1_relaxed} with $\AOSC(\pi^\dagger)=0.4$ and setting the threshold $\kappa$ empirically so that we have $\AOSC(\widehat{\pi})=0.5$. The WADD and ARL for each of the stopping time are computed based on $10^5$ sequences of size $10^6$ generated according to the signal model. It can be seen from \figref{fig:ARLWADDc4} that the GDECuSum and the GLR CuSum with the $2$-threshold stopping time has similar performance. Furthermore, we observe that the GLR Cusum with 2-threshold stopping time $\tau_{\widehat{\pi}}$ has an asymptotic $\ARL$-$\WADD$ trade-off rate that is similar to the GCuSum with the static sampling policy that maximizes the $\ARL$-$\WADD$ trade-off rate.  
		
\begin{figure}[!htb]
		\centering
		\centerline{\includegraphics[width=\columnwidth]{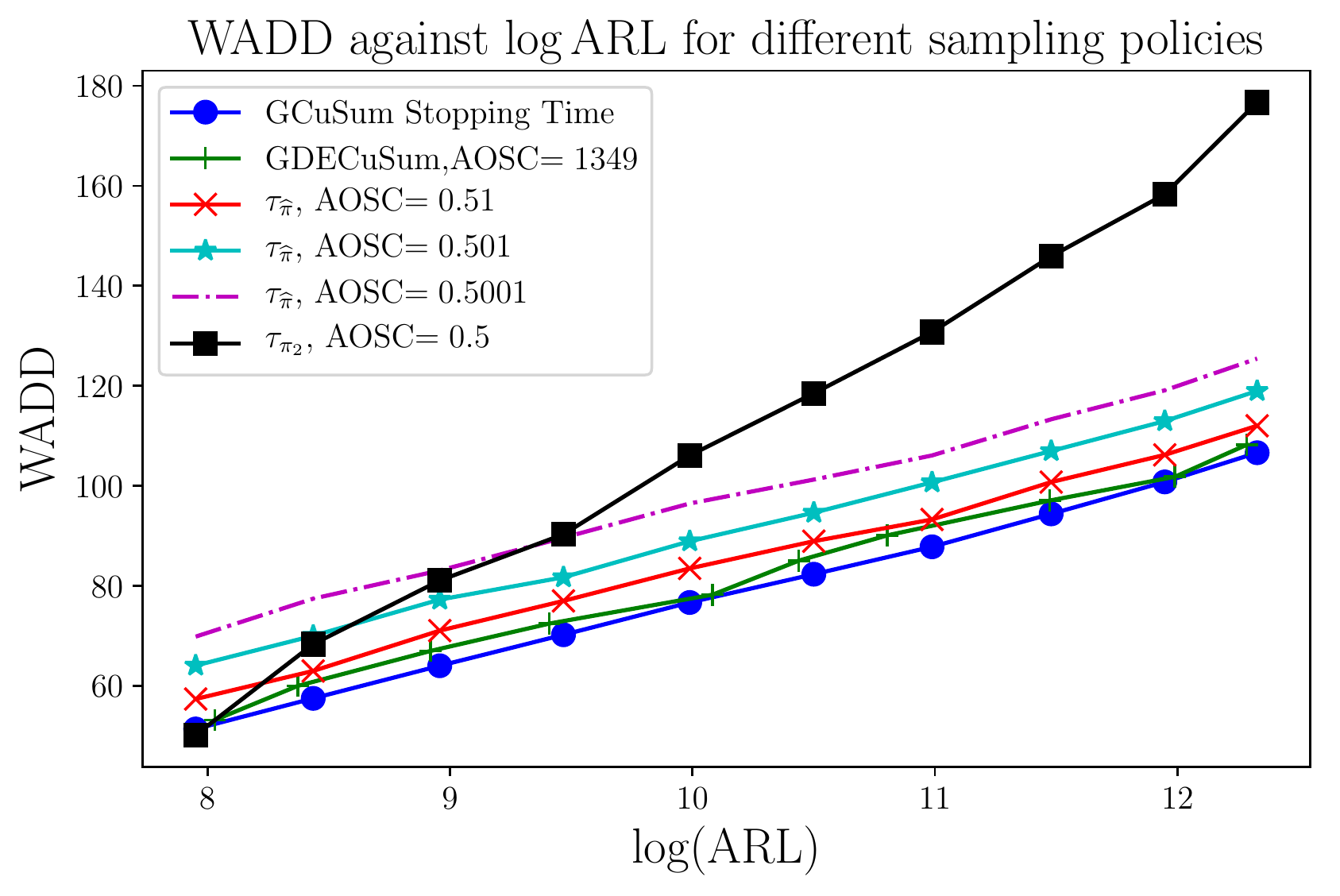}}
		\caption{Comparison of the ARL-WADD trade-off rate for the GDECuSum and the GLR CuSum with the 2-threshold stopping time when $\calC=\calC^5$. }
		\label{fig:ARLWADDc5}
	\end{figure}

		In the second set of simulations, we use another switching cost matrix to illustrate a key difference between our proposed method and the GDECuSum. For gas sensors, the recovery time \cite{lahlalia2019improved,masteghin2020heating} is the time taken for the sensor to reset after taking a measurement. The recovery time for some sensors can be rather long and in some cases, the recovery time may exceed one duty cycle of the sensor. Thus, the sensor is unable to make a reliable observation if it has already made an observation in the previous duty cycle. We can model this phenomenon as having a high action-switching cost when the sensor is switched on consecutively. For the next set of simulations, we use the observation-switching cost matrix 
	$
			\calC^5=\begin{pmatrix}
				10^6 & 0\\
				1&	0
			\end{pmatrix}.
		$  To estimate the $\AOSC$ of the GDECuSum stopping time, we generate $100$ test sequences of size $10^5$ under the pre-change regime. Using these test sequences, the empirical $\AOSC$ of the GDECuSum stopping time is $1349$. This means that the GDECuSum frequently makes consecutive observations which is undesirable. We solve \cref{eqn:lp_policy_design_problem_m=1_relaxed} with $\AOSC=0.5$ and $\epsilon=10^{-10}$ to obtain the $\ARL-\WADD$ trade-off of $0.04$ for our proposed GLR CuSum stopping time with an open-loop policy of window size $1$. Here, it can be seen that our proposed GLR CuSum yields a feasible solution with low $\AOSC$ while the GDECuSum is unable to keep the $\AOSC$ low as the GDECuSum does not take action-switching costs into account. In \figref{fig:ARLWADDc5}, we compare the GCuSum\cite{lorden71} stopping time, GDECuSum\cite{banerjee15} with our proposed GLR CuSum with the 2-threshold stopping time $\tau_{\widehat{\pi}}$ with observation-dependent sampling policy $\widehat{\pi}$. The 2-threshold stopping time is determined by letting $\pi^\dagger$ to be the optimal solution of \cref{eqn:lp_policy_design_problem_m=1_relaxed} with $\AOSC(\pi^\dagger)=0.5$ and setting the threshold $\kappa$ empirically so that we have $\AOSC(\widehat{\pi})=0.51,0.501,0.5001$ respectively. The WADD and ARL for each of the stopping times are computed based on $10^5$ sequences of size $10^6$ generated according to the signal model. It can be seen from \figref{fig:ARLWADDc5} that the GDECuSum and the GCuSum with the static sampling policy that maximizes the $\ARL$-$\WADD$ trade-off rate has similar performance. However, it should be noted that the $\AOSC$ of the GDECuSum is $1349$ as it is unable to take action-switching costs into account. On the other hand, we observe that our proposed GLR Cusum with 2-threshold stopping time $\tau_{\widehat{\pi}}$ has an asymptotic $\ARL$-$\WADD$ trade-off rate that is similar to the GCuSum for $\AOSC=0.51,0.501,0.5001$.  

	\section{Conclusion and Future Work}\label{sec:conclusion}
	In this paper, we discussed the problem of QCD while taking sampling and switching costs into consideration. We formulated the QCD problem with an additional $\AOSC$ constraint. Asymptotically optimal stopping times were proposed and the design of optimal open-loop policies were formulated as quadratic optimization problems. We showed that open-loop policies of window size $W>1$ can be reduced to an open-loop policy of window size $W=1$ while maintaining the same asymptotic $\ARL$-$\WADD$ trade-off and $\AOSC$. Thus, it is sufficient to solve the policy design problem for $W=0,1$. We applied the IRM algorithm to the policy design problem to obtain locally optimal solutions. For cases with additional assumptions on the observation-switching cost matrix $\calC$, globally optimal solutions can be obtained. The methods developed are for the case when the window $W$ is finite. The results regarding the structure of asymptotically optimal stopping times such as \cref{prop:M=1_suffices}, do not hold in general when $W=\infty$ and would be an interesting direction for future work.
	We have also proposed a 2-threshold stopping time $\tau_{\widehat{\pi}}$ for the case when observation-dependent control of the sampling is possible. We showed that with careful construction, the observation-dependent sampling policy $\widehat{\pi}$ together with the GLR CuSum stopping time satisfies the $\AOSC$ and $\ARL$ constraint. Experiment results indicate that the $\ARL$-$\WADD$ trade-off rate of our proposed stopping time is similar to the GCuSum. The derivation of the theoretical bound for the $\WADD$ would be an interesting direction for future work.

	\appendices

	\section{Proof of \cref{prop:information_number_single_recurrence_as_convergence}}\label[appendix]{sec:AppLem2}
	Since $\supp(q)$ lies in the recurrence class $\calR$,	for any $t\geq \nu$ and $\beta\in\calA^W$ such that $q(\beta)>0$, we have $\beta\in \calR$ and 
	\begin{align}
		&\frac{1}{t}\Lambda_m(\nu,\nu+t-1)\nn
		&=\sum_{\alpha\in \calA^W}\frac{1}{t}\sum_{\left\{\substack{i \text{ s.t. }\nu \leq i\leq \nu+t-1, \\  A^{{i-W+1}:i}=\alpha}\right\}}\log\frac{p_m(Y_i|\alpha[W])}{p_0(Y_i|\alpha[W]))}\nn
		&=\sum_{\alpha\in \calR}\frac{N_t(\alpha;\beta)}{t}\frac{1}{N_t(\alpha;\beta)}\sum_{\left\{\substack{i \text{ s.t. }\nu \leq i\leq \nu+t-1, \nn A^{{i-W+1}:i}=\alpha}\right\}}\log\frac{p_m(Y_t|\alpha[W])}{p_0(Y_t|\alpha[W])},\label{aveLambda}
	\end{align}
	where the last equality follows because for any $\alpha\notin \calR$, $\left\{\substack{i \text{ s.t. }\nu \leq i\leq \nu+t-1, \\  A^{{i-W+1}:i}=\alpha}\right\}=\emptyset$. We have
	\begin{align}
		\frac{N_t(\alpha;\beta)}{t}\asure{m} \overline{q}(\alpha)\quad\text{as $t\to \infty$,}
	\end{align}
	and
	\begin{align}
		&\frac{1}{N_t(\alpha;\beta)}\sum_{\left\{\substack{i \text{ s.t. }\nu \leq i\leq \nu+ t-1,\nonumber \\  A^{{i-W+1}:i}=\alpha}\right\}}\log\frac{p_m(Y_t|\alpha[W])}{p_0(Y_t|\alpha[W])}\\
		&\asure{m} \KLD{p_m(\,\cdot\mid\alpha[W])}{p_0(\,\cdot\mid\alpha[W])}\quad\text{as $t\to \infty$.}
	\end{align}
	Thus from \cref{aveLambda}, we obtain 
	\begin{align}
		&\frac{1}{t}\Lambda_m(\nu,\nu+t-1)\nonumber\\
		&\asure{m}\sum_{\alpha\in \calR} \overline{q}(\alpha) \KLD{p_m(\,\cdot\mid\alpha[W])}{p_0(\,\cdot\mid\alpha[W])}=I_{m,\pi}.
	\end{align}
	The proof is now complete.
	
	\section{Proof of \cref{prop:information_number_single_recurrence_lower_bounds_LLR}}\label[appendix]{sec:AppLem3}
	
	Suppose $\supp(\overline{q}) \subseteq \calR$ for a recurrence class $\calR$, then $\supp(q) \subseteq \calR\cup\calU$, where $\calU$ is a set of transient states such that the first-passage probability of entering $\calR$ from each $\beta\in\calU$ is one. Then, from \cref{prop:information_number_single_recurrence_as_convergence}, \cref{eqn:information_number_convergence_in_prob} follows. Let $0\leq j<t $ and $E_j(\nu)$ be the event $\left\{t^{-1}\max_{0\leq j<t}\Lambda_m(\nu,\nu+j) >(1+\epsilon)I_{m,\pi}\right\}$ and we have
	\begin{align}
		&\sup_{0\leq \nu<\infty} \esssup \mathbb{P}_{\nu,m}\left(E_j(\nu)\ \bigg| \ {A^{1:\nu-1},Y^{1:\nu-1}}\right)\nn
		&=\sup_{0\leq \nu<\infty} \esssup \mathbb{P}_{\nu,m}\left(E_j(\nu)\ \bigg| \ {A^{\nu-W:\nu-1}}\right)\label{eqn:prop_3_proof_eqn_2}\\
		&=\sup_{0\leq \nu<\infty} \max_{\substack{\alpha\in\calA^W \\ q(\alpha)>0}} \mathbb{P}_{\nu,m}\left(E_j(\nu)\ \bigg| \ {A^{\nu-W:\nu-1}=\alpha}\right)\nn
		&=\max_{\substack{\alpha\in\calA^W \\ q(\alpha)>0}}  \mathbb{P}_{W+1,m}\left(E_j(W+1)\ \bigg| \ {A^{1:W}=\alpha}\right)\nn
		&\to 0\quad\text{as $t\to\infty$,} \label{eqn:convergence_to_zero_SLLN}
	\end{align}
	where \cref{eqn:prop_3_proof_eqn_2} is because $\Lambda_m(\nu,\nu+j)$ is independent of $A^{1:\nu-W-1}$ and $Y^{1:\nu-1}$ given $A^{\nu-W:\nu-1}$ for each $0\leq j<t$, and \cref{eqn:convergence_to_zero_SLLN} is because each of the terms within the set $\set{\alpha\in\calA^W}|{q(\alpha)>0}$ that we take maximum over converges to zero due to \cref{prop:information_number_single_recurrence_as_convergence}.
	\section{Proof of \cref{prop:reduction_of_recurrence_classes}}\label[appendix]{sec:AppProp2}
	Without loss of generality, we consider the case where $\overline{q}$ has support in the recurrence classes $\calR_1$ and $\calR_2$. Let $q_1,q_2$ be stationary distributions on $\calA^W$ such that $q_i$ has support only in $\calR_i$ for $i\in\{1,2\}$ and that there exists $0<\lambda< 1$ such that 
	$
	\overline{q}=\lambda q_1+(1-\lambda)q_2.
	$ 
	Let $\pi_1=(q_1,T)$ and $\pi_2=(q_2,T)$. From \cref{lem:AOSC_staionary}, we obtain 
	\begin{align}
		\AOSC(q,T)&=\AOSC(\overline{q},T)\nonumber\\
		&=\lambda\AOSC(q_1,T)+(1-\lambda)\AOSC(q_2,T)\nonumber\\
		&\geq \min\left\{\AOSC(\pi_1),\AOSC(\pi_2)\right\}.
	\end{align}
	Let $E_1$ be the event that $\calR_1$ is visited before $\calR_2$ and $E_2$ be the event that $\calR_2$ is visited before $\calR_1$. We have 
	\begin{align}
		&\WADD(\tau,\pi)\nonumber\\
		&=\sup_{\mathclap{\substack{\nu\geq 1\\1\leq m\leq M}}}\esssup \E{\nu,m}[(\tau-\nu+1)^+]{A^{1:\nu-1},Y^{1:\nu-1}}\nonumber\\
		&=\sup_{\mathclap{\substack{\nu\geq 1\\1\leq m\leq M}}}\max_{i=1,2} \esssup\E{\nu,m}[(\tau-\nu+1)^+]{A^{1:\nu-1},Y^{1:\nu-1},E_i}\nonumber\\
		&=\max_{i=1,2}\sup_{\mathclap{\substack{\nu\geq 1\\1\leq m\leq M}}} \esssup\E{\nu,m}[(\tau-\nu+1)^+]{A^{1:\nu-1},Y^{1:\nu-1},E_i}\nonumber\\
		&=\max_{i=1,2}\sup_{\mathclap{\substack{\nu\geq 1\\1\leq m\leq M}}} \esssup\E{\nu,m}[(\tau-\nu+1)^+]{A^{1:\nu-1},Y^{1:\nu-1},E_i}\nonumber\\
		&=\max_{i=1,2}\WADD(\tau,\pi_i).\nonumber
	\end{align} 
	Setting $\pi'=\pi_j$ where $j=\argmin_{i=1,2} \AOSC(\pi_i)$ completes the proof.

	\section{Proof of \cref{thm:structure_of_asym_opt_solutions}}\label[appendix]{sec:AppThm1}
	We first note that \cref{eqn:QCD_AOSC_formulation} is feasible and an optimal solution exists. Let the open-loop policy $\pi^*=(q,T)$ and stopping-time $\tau^*$ such that $(\tau^*,\pi^*)$ is optimal for \cref{eqn:QCD_AOSC_formulation}. From \cref{prop:reduction_of_recurrence_classes}, we can find an open-loop policy $\pi=(q',T)$ such that $\overline{q'}$ has support in only one recurrence class and 
	\begin{align*}
		\AOSC(\pi) &\leq \AOSC(\pi^*),\\ 
		\WADD(\tau^*,\pi) &\leq \WADD(\tau^*,\pi^*).
	\end{align*}
	From the discussion before \cref{eqn:QCD_formulation_specific_sampling}, the GLR CuSum stopping time $\tau_{\pi}$ is asymptotically optimal for the following problem: 
	\begin{equation}\label{eqn:QCD_formulation_specific_sampling_1}
	\begin{aligned}
	& \min_{\tau} &  & \WADD(\tau,\pi)& \st       &  & \ARL(\tau,\pi)\geq \gamma.
	\end{aligned}
	\end{equation}
	Thus, the GLR CuSum stopping time $\tau_{\pi}$ satisfies
	\begin{align*}
		\WADD(\tau_\pi,\pi)\asymp\WADD(\tau^*,\pi)\leq \WADD(\tau^*,\pi^*),
	\end{align*}
	as $\gamma\to\infty$. Hence, $(\tau_{\pi},\pi)$ is also an asymptotically optimal solution to \cref{eqn:QCD_AOSC_formulation}.
	
	\section{Proof of \cref{prop:M=0_suffices}}\label[appendix]{sec:AppProp4}
	If $\pi$ is equivalent to an open-loop policy of window size $0$, then $\pi_0=\pi$ proves the proposition. If $\pi=(q,T)$ is an open-loop policy of window size $W>0$, by \cref{thm:structure_of_asym_opt_solutions}, it suffices to consider the case where the support of  $\overline{q}$ is a subset of a recurrence class of $T$ . Let $\pi_0=(q_0,T_0)$ such that $q_0=\overline{q}^0$ and $T_0$ be the probability transition matrix with rows equal to $q_0$ where \begin{align*}
		\overline{q}^0(a)&=\sum_{{\alpha\in\calA^W : \alpha[W]=a}}\overline{q}(\alpha)
	\end{align*}
	for any $a\in\calA$. For any $m\in\{1,\ldots,M\}$, we have 
	\begin{align}
		I_{m,\pi}&=\sum_{\alpha\in \calA^W}\overline{q}(\alpha) \KLD{p_m(\,\cdot\mid\alpha[W])}{p_0(\,\cdot\mid\alpha[W])}\nonumber\\
		&=\sum_{a\in \calA}\sum_{\left\{\substack{\alpha\in \calA^W :\\ \alpha[W]=a}\right\}}\overline{q}(\alpha) \KLD{p_m(\,\cdot\mid\alpha[W])}{p_0(\,\cdot\mid\alpha[W])}\nonumber\\
		&=\sum_{a\in \calA} \KLD{p_m(\,\cdot\mid a)}{p_0(\,\cdot\mid a)}\sum_{\left\{\substack{\alpha\in \calA^W:\\ \alpha[W]=a}\right\}}\overline{q}(\alpha)\nonumber\\
		&=\sum_{a\in \calA}q_0(a)\KLD{p_m(\,\cdot\mid a)}{p_0(\,\cdot\mid a)}=I_{m,\pi_0},\label{eqn:due_to_expected_visit_times_for_degenerate_MC}
	\end{align}
	where \cref{eqn:due_to_expected_visit_times_for_degenerate_MC} is due to the fact that $\overline{q_0}=q_0$ as $\pi_0$ is an open-loop policy of window size $0$. By \cref{prop:Information_number_WADD}, we obtain $\WADD(\tau_{\pi},\pi)\asymp\WADD(\tau_{\pi_0},\pi_0)$ as $\gamma\to\infty$. The proof is now complete.
	
		\section{Proof of \cref{prop:equivalence_of_optimization_problem}}\label[appendix]{sec:AppProp4.5}
		We first proof claim $(i)$. Let $(q^*,z^*)$ be an optimal solution to \cref{eqn:policy_design_problem_m=0_LP}. It can be easily checked that $q^*$ is a feasible solution for \cref{eqn:policy_design_problem_m=0}. We prove the statement by contradiction. Suppose $q^*$ is not optimal for \cref{eqn:policy_design_problem_m=0}. Then there exists a feasible $q'$ such that 
		\begin{align*}
			&\max_{1\leq m\leq M} \parens*{\sum_{a\in \mathcal{A}} q'(a) \KLD{p_0(\,\cdot\mid a)}{p_m(\,\cdot\mid a)}}^{-1}\nonumber\\
			&\quad< \max_{1\leq m\leq M} \parens*{\sum_{a\in \mathcal{A}} q^*(a) \KLD{p_0(\,\cdot\mid a)}{p_m(\,\cdot\mid a)}}^{-1}.
		\end{align*}
		Let \begin{align*}
			z'=-1/ \max_{1\leq m\leq M} \parens*{\sum_{a\in \mathcal{A}} q'(a) \KLD{p_0(\,\cdot\mid a)}{p_m(\,\cdot\mid a)}}^{-1}.
		\end{align*}
		This gives
		\begin{align*}
			\min_{1\leq m\leq M} \sum_{a\in \mathcal{A}} q'(a) \KLD{p_0(\,\cdot\mid a)}{p_m(\,\cdot\mid a)}+z'=0.
		\end{align*}
		Thus, $z'$ satisfies $ \sum_{a\in \mathcal{A}} q'(a) \KLD{p_0(\,\cdot\mid a)}{p_m(\,\cdot\mid a)}+z'\geq 0$ for all $m\in\{1,...,M\}$ and $(q',z')$ is a feasible solution with $z'<z^*$. This contradicts the optimality of $(q^*,z^*)$. Hence, $q^*$ is an optimal solution to \cref{eqn:policy_design_problem_m=0}.\\
		
		We now proof claim $(ii)$. Let $q^*$ be an optimal solution to \cref{eqn:policy_design_problem_m=0}. Let \begin{align*}
			z^*=-1/ \max_{1\leq m\leq M} \parens*{\sum_{a\in \mathcal{A}} q'(a) \KLD{p_0(\,\cdot\mid a)}{p_m(\,\cdot\mid a)}}^{-1}.
		\end{align*} Following the manipulation above,  $z^*$ satisfies $ \sum_{a\in \mathcal{A}} z^*(a) \KLD{p_0(\,\cdot\mid a)}{p_m(\,\cdot\mid a)}+z^*\geq 0$ for all $m\in\{1,...,M\}$. Thus, $(q^*,z^*)$ is a feasible solution. We prove the statement by contradiction. Suppose $(q^*,z^*)$ is not optimal for \cref{eqn:policy_design_problem_m=0_LP}. Then there exists a feasible $(q',z')$ such that $z'<z^*$. Since $(q',z')$ is feasible, $ \sum_{a\in \mathcal{A}} q'(a) \KLD{p_0(\,\cdot\mid a)}{p_m(\,\cdot\mid a)}+z'\geq 0$ holds for all $m\in\{1,\ldots,M\}$. Thus, with some algebraic manipulation,
		\begin{align*}
			\max_{1\leq m\leq M} \parens*{\sum_{a\in \mathcal{A}} q'(a) \KLD{p_0(\,\cdot\mid a)}{p_m(\,\cdot\mid a)}}^{-1}\leq \frac{-1}{z'}.
		\end{align*} Since $z'<z^*$, we have $\frac{-1}{z'}<\frac{-1}{z^*}$. By the choice of $z^*$, we have
		\begin{align*}
			&\max_{1\leq m\leq M} \parens*{\sum_{a\in \mathcal{A}} q'(a) \KLD{p_0(\,\cdot\mid a)}{p_m(\,\cdot\mid a)}}^{-1}\nonumber\\
			&\quad\leq \frac{-1}{z'}\nonumber\\
			&\quad< \frac{-1}{z^*}=\max_{1\leq m\leq M} \parens*{\sum_{a\in \mathcal{A}} q^*(a) \KLD{p_0(\,\cdot\mid a)}{p_m(\,\cdot\mid a)}}^{-1}.
		\end{align*} This contradicts the optimality of $q^*$. Hence, $(q^*,z^*)$ is an optimal solution to \cref{eqn:policy_design_problem_m=0_LP}.		
		The proof is now complete.
	
	\section{Proof of \cref{prop:M=1_suffices}}\label[appendix]{sec:AppProp5}
	If $\pi$ is equivalent to an open-loop policy of window size $1$, we let $\pi_1=\pi$ and we have proved the proposition. If $\pi$ is an open-loop policy of window size $W>1$, by \cref{thm:structure_of_asym_opt_solutions}, it suffices to consider the case where the support of  $\overline{{q}}$ is a subset of a recurrence class $\calR$. The $\AOSC$ of $\pi$ can be expressed as 
	\begin{align}
		&\AOSC(\pi)\nonumber\\
		&=\AOSC(\overline{q},T)\nonumber\\
		&=\limsup_{n\to \infty}\frac{1}{n}\E{\infty}[\sum_{t=2}^{n+1}\calC[A_{i-1},A_{i}]]\nonumber\\
		&=\limsup_{n\to \infty}\frac{1}{n}\sum_{\beta\in\calA^W}\overline{q}(\beta)\E{\infty}[\sum_{t=2}^{n+1}\calC[A_{i-1},A_{i}]]{A^{1:W}=\beta}\nonumber\\
		&=\limsup_{n\to \infty}\frac{1}{n}\sum_{\beta\in\calA^W}\overline{q}(\beta)\E{\infty}[\sum_{\alpha\in\calA^W}\calC_\alpha N_{n+1}(\alpha;\beta)]{A^{1:W}=\beta}\nonumber\\
		&=\sum_{\beta\in\calA^W}\sum_{\alpha\in\calA^W}\overline{q}(\beta)\calC_\alpha\overline{q}(\alpha)=\sum_{\alpha\in\calA^W}\calC_\alpha\overline{q}(\alpha).\nonumber
	\end{align}
	Let $\pi_2=(q_2,T_2)$ be an open-loop policy of window size $2$ and it can be similarly shown that
	$
	\AOSC(\pi_2)=\sum_{\alpha\in\calA^2}\calC_\alpha\overline{q_2}(\alpha).
	$
	Thus, it can be seen that the $\AOSC(\pi_2)$ depends only on $\overline{q_2}$. Let $q_2=\overline{q}^1$ be the projection of $\overline{q}$ onto $\calA^2$ where 
	\begin{align*}
		\overline{q}^1(a,b)&=\sum_{\substack{\alpha\in\calA^W : \\ \alpha[W-1]=a,\alpha[W]=b}}\overline{q}(\alpha),
	\end{align*}
	for any $a,b\in\calA$ and $T_2$ be defined as
	\begin{align*}
		T_2[(a,b),(c,d)]=\begin{cases}
			\frac{\overline{q}^1(c,d)}{\sum_{d\in\calA}\overline{q}^1(c,d)}\quad\text{if $c=b$,}\\
			0\quad\text{otherwise.}
		\end{cases}
	\end{align*} 
	For any $(a,b)\in\calA^2$, we have
	\begin{align*}
		\sum_{(c,d)\in\calA^2}	T_2[(a,b),(c,d)]&=\sum_{d\in\calA}\frac{\overline{q}^1(b,d)}{\sum_{d_1\in\calA}\overline{q}^1(b,d_1)}\nonumber\\
		&=\frac{\sum_{d\in\calA}\overline{q}^1(b,d)}{\sum_{d_1\in\calA}\overline{q}^1(b,d_1)}=1,
	\end{align*}
	and hence, $T_2$ is a probability transition matrix. Since the support of $\overline{q}$ is a subset of a single recurrence class of $T$, the support of $q_2$ also lies in a single recurrence class of $T_2$. Next, we show that $q_2$ is a stationary distribution of $T_2$:
	\begin{align*}
		T_2q_2[(c,d)]&=\sum_{(a,b)\in\calA^2}T_2[(a,b),(c,d)]\overline{q}^1(a,b)\nonumber\\
		&=\sum_{a\in\calA}T_2[(a,c),(c,d)]\overline{q}^1(a,c)\nonumber\\
		&=\sum_{a\in\calA}\frac{\overline{q}^1(c,d)}{\sum_{d\in\calA}\overline{q}^1(c,d)}\overline{q}^1(a,c)\nonumber\\
		&=\frac{\overline{q}^1(c,d)}{\sum_{d\in\calA}\overline{q}^1(c,d)}\sum_{a\in\calA}\overline{q}^1(a,c)\nonumber\\
		&=\overline{q}^1(c,d)=q_2(c,d).
	\end{align*}
	Therefore, $q_2$ is a stationary distribution of $T_2$ which has support contained in one recurrence class. Hence, $\overline{q_2}=q_2$ and $\AOSC(\pi)=\AOSC(\pi_2)$. Using similar arguments from the proof of \cref{prop:M=0_suffices}, we have $I_{m,\pi}=I_{m,\pi_1}$ for $m=\{1,\ldots,M\}$ and $\WADD(\tau_{\pi},\pi)\asymp\WADD(\tau_{\pi_2},\pi_2)$ as $\gamma\to\infty$.
	
	Furthermore, a quick computation shows that $\pi_2=(q_2,T_2)$ is equivalent to the open-loop policy $\pi_1=(q_1,T_1)$ of window size $1$ with 
	$
	q_1(a)=\sum_{b \in\calA}q_2(a,b)$ and $T_1[a,b]=\frac{q_2(a,b)}{\sum_{b \in\calA}q_2(a,b)}.
	$
	
	Thus, we have $\AOSC(\pi)=\AOSC(\pi_1)$ and $\WADD(\tau_{\pi},\pi)\asymp\WADD(\tau_{\pi_1},\pi_1)$ as $\gamma\to\infty$. The proof is complete.

		\section{Counterexample to \cref{prop:M=1_suffices} for $W=\infty$}\label[appendix]{sec:AppExp1}
		
		We consider the problem with two actions $\calA=\{1,2\}$ and a $2\times2$ observation-switching costs matrix $\mathcal{C}=\mathcal{C}^*$ as follows:
		\begin{align*}
			\mathcal{C}^*[i,j]=\begin{cases}
				0 \quad\text{if $i=j$}\\
				1 \quad\text{if $i\neq j$},
			\end{cases}
		\end{align*}
		Let $M=2$. The pre- and post-change distributions are chosen with the following pmfs:
		\begin{align*}
			p_0=\left[\frac{1}{3},\frac{1}{3},\frac{1}{3}\right],\ p_1=\left[\frac{1}{3},\frac{1}{6},\frac{1}{2}\right],\ p_2=\left[\frac{1}{2},\frac{1}{6},\frac{1}{3}\right].
		\end{align*}
		
		We define the actions so that given the actions, pre- and post-change distributions are discrete distributions with the following pmfs:
		\begin{align*}
			&\P{\nu,m}(Y_t=1){A=1}=p_m[1]+p_m[2],\\
			&\P{\nu,m}(Y_t=2){A=1}=p_m[3],\\
			&\P{\nu,m}(Y_t=1){A=2}=p_m[1],\\
			& \P{\nu,m}(Y_t=2){A=2}=p_m[2]+p_m[3],
		\end{align*}for $t\geq\nu$ and $m\in\{1,2\}$,
		\begin{align*}
		&\P{\nu,m}(Y_t=1){A=1}=p_0[1]+p_0[2],\\
		&\P{\nu,m}(Y_t=2){A=1}=p_0[3],\\
		&\P{\nu,m}(Y_t=1){A=2}=p_0[1],\\
		& \P{\nu,m}(Y_t=2){A=2}=p_0[2]+p_0[3],
		\end{align*}for $t<\nu$, and
		\begin{align*}
		&\P{\infty}(Y_t=1){A=1}=p_0[1]+p_0[2],\\
		&\P{\infty}(Y_t=2){A=1}=p_0[3],\\
		&\P{\infty}(Y_t=1){A=2}=p_0[1],\\
		& \P{\infty}(Y_t=2){A=2}=p_0[2]+p_0[3],
	\end{align*}
	for all $t\in\mathbb{N}$.
		To compute the upper-bound of the ARL-WADD trade-off, we consider the QCD with observation-switching costs problem with $\calC=\calC^0=\mathbf{0}$ where $\mathbf{0}$ is the $2\times2$ zero matrix. The ARL-WADD trade-off when $\calC=\calC^0$ is $\frac{1}{4}\left(\log\frac{3}{4}+\log\frac{3}{2}\right)$.
		
		We consider a deterministic policy $\pi_\infty$(i.e., Markov chain with $W=\infty$) and sample the signal using the following sequence of actions:
		\begin{align*}
			\{A_t\}=\{1,2,1,1,2,2,1,1,1,2,2,2,1,1,1,1,2,2,2,2,...\}
		\end{align*}
		where the number of times each action is used to sample the signal increases by 1 before switching. When $\calC=\mathcal{C}^*$, the average number of switches tends to zero as $t\to \infty$ and we have $\text{ASOC}(\pi_\infty)=0$ and 
		\begin{align*}
			\P{\nu,m}(\lim_{t\to\infty}\frac{1}{t}\Lambda_m(\nu,\nu+t)=\frac{1}{4}\left(\log\frac{3}{4}+\log\frac{3}{2}\right))=1,
		\end{align*}for any $\nu\in\mathbb{N}$ and $m\in\{1,2\}$. Thus, the policy $\pi_\infty$ achieves asymptotically optimal $\WADD$-$\ARL$ trade-off with $\text{ASOC}(\pi_\infty)=0$.
		
		However, $\text{ASOC}(\pi_1)>0$ for any policy $\pi_1$ of window size $W=1$ where both actions $1$ and $2$ are used. On the other hand, if only one of the actions is used, then either $I_{1,\pi_1}$ or $I_{2,\pi_1}$ is zero. Hence, $\pi_\infty$ cannot be reduced to a policy of window size $W=1$.

	\AtEndDocument{\refstepcounter{equation}\label{finalthm}}
\nocite{chernoff1959sequential}
\nocite{nitinawarat2013controlled}
	\bibliographystyle{IEEEtran}
	\bibliography{IEEEabrv,StringDefinitions,refs}

\title{Supplementary Material for "Asymptotically Optimal Sampling Policy for Quickest Change Detection with Observation-Switching Cost"}
\author{
	Tze~Siong~Lau and Wee~Peng~Tay,~\IEEEmembership{Senior Member,~IEEE}
	\thanks{This work was supported in part by the Singapore Ministry of Education Academic Research Fund Tier 2 grant MOE2018-T2-2-019 and by A*STAR under its RIE2020 Advanced Manufacturing and Engineering (AME) Industry Alignment Fund – Pre Positioning (IAF-PP) (Grant No. A19D6a0053).}%
	\thanks{T.~S.~Lau and W.~P.~Tay are with the School of Electrical and Electronic Engineering, Nanyang Technological University, Singapore (e-mail: TLAU001@e.ntu.edu.sg, wptay@ntu.edu.sg).}
}
%


%

	\maketitle
	\begin{abstract}
		Supplementary material for the paper "Asymptotically Optimal Sampling Policy for Quickest Change Detection with Observation-Switching Cost".
	\end{abstract}
	%
	\begin{IEEEkeywords}
		Quickest change detection, GLR CuSum, sampling policy, graph sampling
	\end{IEEEkeywords}

	The notations and references in this supplementary material refer to those in the paper "Asymptotically Optimal Sampling Policy for Quickest Change Detection with Observation-Switching Cost". 
	\section{Proof of Lemma 1}\label[appendix]{sec:AppLem1}
	
	Each time a state $\alpha$ is visited, a switching cost $\calC_\alpha$ is incurred. Thus, the \gls{AOSC} can be written as 
	\begin{align}\label{eqn:AOSC_pi1}
		\AOSC(\pi_1)&=\lim_{t\to\infty}\sum_{\beta\in\calA^W}\sum_{\alpha\in\calA^W}\E{\infty}[\frac{\calC_\alpha N_t(\alpha;\beta)}{t}q(\beta)]\nonumber\\
		&=\sum_{\beta\in\calA^W}\sum_{\alpha\in\calA^W}q(\beta)\calC_\alpha\lim_{t\to\infty}\E{\infty}[\frac{N_t(\alpha;\beta)}{t}].
	\end{align}
	Similarly, for $\pi_2$, we have 
	\begin{align}\label{eqn:AOSC_pi2}
		\AOSC(\pi_2)&=\sum_{\beta\in\calA^W}\sum_{\alpha\in\calA^W}\overline{q}(\beta)\calC_\alpha\lim_{t\to\infty}\E{\infty}[\frac{N_t(\alpha;\beta)}{t}].
	\end{align}
	Suppose there are $R$ recurrence classes $\calR_1,\ldots,\calR_R$. For any $r\in\{1,\ldots, R\}$ with both $\beta,\delta\in \calR_r$, we have 
	\begin{align}\label{eqn:equal_average_in_same_class}
		\lim_{t\to\infty}\E{\infty}[\frac{N_t(\alpha;\beta)}{t}]=\lim_{t\to\infty}\E{\infty}[\frac{N_t(\alpha;\delta)}{t}].
	\end{align}
	Let $\beta_r\in \calR_r$ for $r=1,\ldots, R$ and for any transient state $\delta$, we have 
	\begin{align}\label{eqn:decomposition_in_transient_class}
		\lim_{t\to\infty}\E{\infty}[\frac{N_t(\alpha;\delta)}{t}]&=\sum_{r=1}^Rf_{\delta,r}\lim_{t\to\infty}\E{\infty}[\frac{N_t(\alpha;\beta_r)}{t}].
	\end{align}
	Let $\{\delta_j\ |\ j=1,\ldots,U\}$ be the transient states. For any $r=1,\ldots,R$, we have 
	\begin{align}\label{eqn:stationary_distribution}
		\sum_{\beta\in \calR_r}\overline{q}(\beta)=\sum_{\beta\in \calR_r}q(\beta)+\sum_{j=1}^U f_{\delta_j,r}q(\delta_j).
	\end{align}
	Combining equations \cref{eqn:AOSC_pi1,eqn:AOSC_pi2,eqn:equal_average_in_same_class,eqn:decomposition_in_transient_class,eqn:stationary_distribution} together, yields
	\begin{align}
			&\AOSC(\pi_2)\nn
			&=\sum_{\beta\in\calA^W}\sum_{\alpha\in\calA^W}\overline{q}(\beta)\calC_\alpha\lim_{t\to\infty}\E{\infty}[\frac{N_t(\alpha;\beta)}{t}]\\
			&=\sum_{r=1}^R\sum_{\beta\in \calR_r}\sum_{\alpha\in\calA^W}\overline{q}(\beta)\calC_\alpha\lim_{t\to\infty}\E{\infty}[\frac{N_t(\alpha;\beta)}{t}]\label{eqn:AOSC_staionary_1}\\
			&=\sum_{r=1}^R\sum_{\beta\in \calR_r}\sum_{\alpha\in \calR_r}\overline{q}(\beta)\calC_\alpha\lim_{t\to\infty}\E{\infty}[\frac{N_t(\alpha;\beta)}{t}]\label{eqn:AOSC_staionary_1_2}\\
			&=\sum_{r=1}^R\sum_{\beta\in \calR_r}\overline{q}(\beta)\sum_{\alpha\in \calR_r}\calC_\alpha\xi_r[\alpha]\label{eqn:AOSC_staionary_1_3}\\
			&=\sum_{r=1}^R\sum_{\beta\in \calR_r}q(\beta)\sum_{\alpha\in \calR_r}\calC_\alpha\xi_r[\alpha]+\sum_{r=1}^R\sum_{j=1}^U f_{\delta_j,r}q(\delta_j)\sum_{\alpha\in \calR_r}\calC_\alpha\xi_r[\alpha]\label{eqn:AOSC_staionary_1_4}\\
			&=\sum_{r=1}^R\sum_{\beta\in \calR_r}\sum_{\alpha\in \calR_r}q(\beta)\calC_\alpha\lim_{t\to\infty}\E{\infty}[\frac{N_t(\alpha;\beta)}{t}]+\sum_{r=1}^R\sum_{j=1}^U \sum_{\alpha\in \calR_r}f_{\delta_j,r}q(\delta_j)\calC_\alpha\xi_r[\alpha]\label{eqn:AOSC_staionary_1_5}
		\end{align}
		where \cref{eqn:AOSC_staionary_1} is due to $\xi_r[\alpha]=0$ for $\alpha\notin \calR_r$, \cref{eqn:AOSC_staionary_1_4} is obtained by applying \cref{eqn:stationary_distribution} to \cref{eqn:AOSC_staionary_1_3} and \cref{eqn:AOSC_staionary_1_5} is obtained by applying \cref{eqn:equal_average_in_same_class} to \cref{eqn:AOSC_staionary_1_4}.  We also have 
		\begin{align}
			&\AOSC(\pi_1) \nn
			&=\sum_{\beta\in\calA^W}\sum_{\alpha\in\calA^W}q(\beta)\calC_\alpha\lim_{t\to\infty}\E{\infty}[\frac{N_t(\alpha;\beta)}{t}]\label{eqn:AOSC_staionary_2_1}\\
			&=\sum_{r=1}^R\sum_{\beta\in \calR_r}\sum_{\alpha\in\calA^W}q(\beta)\calC_\alpha\lim_{t\to\infty}\E{\infty}[\frac{N_t(\alpha;\beta)}{t}]+\sum_{j=1}^U\sum_{\alpha\in\calA^W}q(\delta_j)\calC_\alpha\lim_{t\to\infty}\E{\infty}[\frac{N_t(\alpha;\delta_j)}{t}]\label{eqn:AOSC_staionary_2_2}\\
			&=\sum_{r=1}^R\sum_{\beta\in \calR_r}\sum_{\alpha\in\calA^W}q(\beta)\calC_\alpha\lim_{t\to\infty}\E{\infty}[\frac{N_t(\alpha;\beta)}{t}]+\sum_{r=1}^R\sum_{j=1}^U \sum_{\alpha\in\calA^W}f_{\delta_j,r}q(\delta_j)\calC_\alpha\xi_r[\alpha]\label{eqn:AOSC_staionary_2_4}\\
			&=\sum_{r=1}^R\sum_{\beta\in \calR_r}\sum_{\alpha\in\calA^W}q(\delta_j)\calC_\alpha\lim_{t\to\infty}\E{\infty}[\frac{N_t(\alpha;\beta)}{t}]+\sum_{r=1}^R\sum_{j=1}^U \sum_{\alpha\in \calR_r}f_{\delta_j,r}q(\delta_j)\calC_\alpha\xi_r[\alpha],\label{eqn:AOSC_staionary_2}
		\end{align}
		where \cref{eqn:AOSC_staionary_2_2} is obtained by partitioning the actions in $\calA^W$ into recurrence classes $\calR_1,\ldots\calR_R$ and transient states, \cref{eqn:AOSC_staionary_2_4} is obtained by applying \cref{eqn:decomposition_in_transient_class} to \cref{eqn:AOSC_staionary_2_2} and, \cref{eqn:AOSC_staionary_2} is due to $\xi_r[\alpha]=0$ for $\alpha\notin \calR_r$.
	Thus, we have $\AOSC(\pi_1)=\AOSC(\pi_2)$ and the proof is complete.
	\bibliographystyle{IEEEtran}
	\bibliography{IEEEabrv,StringDefinitions,refs}
\end{document}